\documentclass[a4,11pt]{amsart}
\usepackage{amssymb,amsmath,amsthm,txfonts}
\usepackage{cite}

\usepackage{color}

\newtheorem{thm}{Theorem}[section]
\newtheorem{lem}[thm]{Lemma}

\newtheorem{prop}[thm]{Proposition}
\theoremstyle{definition}

\theoremstyle{remark}

 \makeatletter
    
    \@addtoreset{equation}{section}
  \makeatother

\makeatletter
      \def\@makefnmark{%
         \leavevmode
            \raise.9ex\hbox{\check@mathfonts
                \fontsize\sf@size\z@\normalfont%
                            \@thefnmark}%
       }
      \makeatother

\newcommand{\dd}{\textrm{d}}

\begin{document}

\title[]{Existence of vortex rings in Beltrami flows}
\author[]{Ken Abe}
\date{}
\address[K. ABE]{Department of Mathematics, Graduate School of Science, Osaka City University, 3-3-138 Sugimoto, Sumiyoshi-ku Osaka, 558-8585, Japan}
\email{kabe@sci.osaka-cu.ac.jp}

\subjclass[2010]{35Q31, 35Q35}
\keywords{Vortex rings, Beltrami fields, Grad-Shafranov equation}
\date{\today}

\maketitle

\begin{abstract}
We construct traveling wave solutions to the 3d Euler equations by axisymmetric Beltrami fields with a non-constant proportionality factor. They form a vortex ring with nested invariant tori consisting of level sets of the proportionality factor. 
\end{abstract}

\vspace{5pt}

\section{Introduction}

\vspace{5pt}

We consider the 3d Euler equations

\begin{equation*}
\begin{aligned}
 v_t+v\cdot \nabla v+\nabla q=0,\quad \nabla\cdot v=0, \quad x\in \mathbb{R}^{3},\ t>0,
\end{aligned}
\tag{1.1}
\end{equation*}\\ 
where the vector field $v(x,t)$ denotes the velocity of a fluid and the scaler field $q(x,t)$ denotes the pressure. The subscript $v_t$ denotes the partial derivative for $t$ and $\nabla=\nabla_x$ denotes the gradient for $x$. We study traveling wave solutions to (1.1) among classical and axisymmetric solutions.

Traveling waves appear in the study of shallow water equations. A classical model of the shallow water is the KdV equation \cite{KdV}

\begin{align*}
v_t+6vv_x+v_{xxx}=0,\quad x\in  \mathbb{R},\ t>0,    
\end{align*}\\
where the scaler field $v(x,t)$ denotes the hight of a wave above the flat bottom. This equation is a dispersive equation of third order. It is known that the KdV equation is globally well-posed for smooth data, e.g. \cite{KPV93}, and does not exhibit a wave breaking, i.e. a finite time blow-up with bounded $v$ and  unbounded $v_{x}$ \cite[p.476]{Wh}. The KdV equation admits a traveling wave solution $v(x,t)=u(x-ct)$ with the speed $c>0$, called a soliton.

A dispersionless model is the Camassa-Holm equation \cite{CH93}

\begin{align*}
v_t+vv_x+q_x=0,\quad q=\frac{1}{2}e^{-|x|}*\left(v^{2}+\frac{v_x^{2}}{2}\right),\quad x\in \mathbb{R}, \ t>0,   
\end{align*}\\
where the symbol $*$ denotes the convolution. This equation is derived from the 2d Euler equations with gravity and a free surface. In terms of $\zeta=v-v_{xx}$, the Camassa-Holm equation can be written as 

\begin{align*}
\zeta_t+v\zeta_x+2\zeta v_x=0,\quad x\in \mathbb{R}, \ t>0.  
\end{align*}\\
Global existence and a finite time blow-up of solutions to this equation are proved by Constantin and Escher. It is known that the Camassa-Holm equation admits a global-in-time solution for smooth data with non-negative $\zeta_0(x)=\zeta(x,0)$  \cite{CE98a} and exhibits a wave breaking for some $\zeta_0$ changing its sign \cite{CE98b}. A traveling wave is called a peakon \cite{CH93}.

A 2-component generalization \cite{SR04} is 

\begin{align*}
\zeta_t+v\zeta_x+2\zeta v_x+ \rho\rho_x=0,\quad 
\rho_t+(v\rho)_x=0,\quad x\in \mathbb{R},\ t>0.  
\end{align*}\\
The unknown $\rho$ is related with horizontal variation of a free surface. This equation is globally well-posed for small data and exhibits a wave breaking for some large data \cite{CI08}. A traveling wave also exists. Both KdV and (2-component) Camassa-Holm equations are infinite-dimensional bi-Hamiltonian systems \cite[p.438]{Olver}, \cite{CI08}.

The 3d Euler equations can be written in a similar form. In the axisymmetric setting, it is recasted as 

\begin{equation*}
\begin{aligned}
\zeta_t +v\cdot \nabla_x \zeta-\frac{1}{r^{4}}\rho\rho_z=0,\quad 
\rho_t +\nabla_x\cdot ( v\rho)=0, \quad z\in \mathbb{R},\ r>0,\ t>0,
\end{aligned}
\tag{1.2}
\end{equation*}\\ 
where $\zeta(z,r,t)$ and $\rho(z,r,t)$ are related with the $\theta$-component of the vorticity $\nabla \times v$ and the velocity $v$ in the cylindrical coordinate $(r,\theta,z)$. The $\theta$-component of the velocity is called \textit{swirl}. The problem without swirl $\rho\equiv 0$ has a similarity as the 2d Euler equations and is globally well-posed for smooth data \cite{UI}. 

On the other hand, global existence of solutions with swirl $\rho\nequiv 0$ is unknown due to lack of an a priori estimate. Some numerical works detected appearance of singularities though definitive results seem unknown, e.g.  \cite[p.185]{Mab}. It is known \cite{BKM} that a finite time blow-up of solutions is equivalent to unboundedness of vorticity. This unboundedness is prevented if vorticity is smoothly directed \cite{CFM}. 

The equations (1.2) has a similarity as the 2d Bousinessq equations for which $r^{-4}$ is replaced with a constant \cite{Mab}. Besides the axis $r=0$, the interior region $r>0$ is also a candidate of a blow-up point. The numerical simulation \cite{LuoHou14} detected appearance of a singularity for an axisymmetric Euler flow on the boundary of a cylinder. Its 1d model exhibits a finite time blow-up \cite{CHKLSY}. The equation (1.2) is also a Hamiltonian system \cite[p.40]{Ben84}, cf. \cite[p.444]{Olver}. \\

Traveling waves of the Euler equations are related with translating vortex motion dating back to a pioneering work of Helmholtz \cite{Hel58}. They form

\begin{equation*}
\begin{aligned}
v(x,t)=u(x+u_{\infty}t)-u_{\infty},  \quad
q(x,t)=p(x+u_{\infty}t),
\end{aligned}
\end{equation*}\\ 
where $v$ vanishes at space infinity with a constant $u_{\infty}$. By substituting this into the Euler equation, one sees that the profile $(u,p)$ solves the steady equations

\begin{equation*}
\begin{aligned}
(\nabla\times u)\times u+\nabla \Pi=0 ,\quad \nabla\cdot u&=0, \hspace{17pt} \textrm{in}\   \mathbb{R}^{3},\\
u&\to u_\infty\quad \textrm{as}\ |x|\to\infty,
\end{aligned}
\tag{1.3}
\end{equation*}\\ 
with the Bernoulli pressure $\Pi=p+|u|^{2}/2$. This equation can describe a broad class of translating vortices. Because of the 3d nature of the domain, a vortex can be supported on a knotted and linked solid torus. Existence of such a vortex is a conjecture of Kelvin \cite[p.264]{Ricca}.

Besides a support of a vortex, integral curves of the velocity (stream lines) and the vorticity (vortex lines) can be also knotted and linked. If the Bernoulli pressure $\Pi$ is not constant, it acts as a first integral of stream lines and vortex lines. They lie on level sets of $\Pi$ (Bernoulli surfaces). The Bernoulli surfaces are known to be nested tori or surfaces diffeomorphic to cylinders by the structure theorem of Arnold \cite[p.69]{AK98}, though his theorem is stated for steady Euler flows in a bounded domain with analytic velocity. All stream lines on each torus are \textit{closed} or \textit{quasi-periodic}, i.e. dense on the torus. Such a torus is called an \textit{invariant torus} in terms of mechanics \cite[p.272]{Arnold}.

If the Bernoulli pressure is constant, the velocity and the vorticity are collinear, i.e. $u\parallel \nabla\times  u$. Such a vector field is called a \textit{Beltrami field}. With a proportionality factor $f$, the steady Euler equations (1.3) are reduced to

\begin{equation*}
\begin{aligned}
\nabla\times u =f u ,\quad \nabla\cdot u&=0, \hspace{17pt} \textrm{in}\   \mathbb{R}^{3},\\
u&\to u_\infty\quad \textrm{as}\ |x|\to\infty.
\end{aligned}
\tag{1.4}
\end{equation*}\\
If the factor $f$ is not constant, it is a first integral of vortex lines (stream lines) and plays an alternative role of the Bernoulli pressure \cite[p.71]{AK98}. 

An exceptional case is when the factor $f$ is constant. A Beltrami field with a constant factor is called a \textit{strong} Beltrami field. Due to absence of the first integrals, vortex lines of a strong Beltrami field have topological freedom and can be chaotic, e.g. ABC flows \cite{AK98}. It is known \cite{EP12} that for any linked curves, there exists a strong Beltrami field (with $u_{\infty}= 0$) having a vortex line diffeomorphic to the curve. Kelvin's conjecture is revisited by Enciso and Peralta-Salas \cite{EP15} by constructing a strong Beltrami field with a \textit{thin} knotted and linked \textit{vortex tube}. A vortex tube is a union of vortex lines embedded to a solid torus. The strong Beltrami field \cite{EP15} has knotted and linked vortex lines.

A vortex of a strong Beltrami field can not be compactly supported since a strong Beltrami field is an eigenfunction of the Laplace operator $-\Delta u=f^{2}u$, which is trivial if $u\in L^{2}(\mathbb{R}^{3})$ by the Fourier transform. The strong Beltrami fields \cite{EP12}, \cite{EP15} decay by the sharp order $u=O(|x|^{-1})$ as $|x|\to\infty$. It is known \cite{Na14} that Beltrami fields must be trivial if $u\in L^{q}(\mathbb{R}^{3})$ for $q\in [2,3]$ or $u=o(|x|^{-1})$ as $|x|\to\infty$.\\

\vspace{-1pt}
In the axisymmetric setting, a vortex of a steady Euler flow can be supported on an unknotted and unlinked solid torus, usually referred to as a \textit{vortex ring}. A vortex ring is a vortex tube with compactly supported vorticity. Due to the axisymmetry, circulation of the velocity acts as an additional first integral of stream lines. We denote by $2\pi \Gamma$ circulation of the velocity along a circle around the symmetric axis. There are 3 particular types of vortex rings: \\ 

\noindent
\textit{Type I}: $\Gamma\equiv 0$. This is an axisymmetric flow without swirl. The velocity and the vorticity are perpendicular as 2d, i.e. $u\perp \nabla \times u$. Vortex lines are circles, unknotted and unlinked. Existence of this type vortex ring is proved by Fraenkel \cite{Fra70}, \cite{Fra72}.\\

\noindent
\textit{Type II}: $\Pi\equiv \textrm{const}$. This is a Beltrami flow. The factor $f$ is related with the circulation. Vortex lines (stream lines) are \textit{knotted} and \textit{linked}. An explicit solution of this type vortex ring with a non-constant factor is founded by Chandrasekhar \cite{Chandra}. \\

\noindent
\textit{Type III}: $u_{\infty}= 0$. This is a steady flow rest at infinity. The velocity can be also compactly supported besides the vorticity. An explicit solution is found by Prendergast \cite{Prend}.\\

Hicks \cite{Hicks99} and Moffatt \cite{Moffatt} found an explicit solution to the steady Euler equations (1.3), called the \textit{Hicks-Moffatt vortex}, which is a 3 parameter family of solutions with a vortex supported on a \textit{ball} including the above 3 types as particular cases, see Section 2. One of them agrees with the Beltrami field of Chandrasekhar \cite{Chandra}. Besides it, a Beltrami field with a non-constant factor is constructed by Turkington \cite[p.68]{Tu89} by a variational principle.

In contrast to many strong Beltrami fields, existence of a Beltrami field with a non-constant factor seems unknown besides \cite{Chandra}, \cite{Tu89}. Enciso and Peralta-Salas \cite{EP16} proved non-existence of Beltrami fields if $f\in C^{2+\mu}(\mathbb{R}^{3})$ for some $\mu\in (0,1)$ and a level set $f^{-1}(k)$ is diffeomorphic to a \textit{sphere} for some $k\in \mathbb{R}$, e.g. $f$ is radial or has local extrema. The non-existence result might imply that realization of Kelvin's conjecture by Beltrami fields with a non-constant factor is subtle. The Beltrami fields \cite{Chandra}, \cite{Tu89} have discontinuous factors $f$ whose level set $f^{-1}(k)$ is a \textit{ball} or a \textit{solid torus} for some $k$ and its complement or an empty set for other $k$, see Section 2.

The aim of this paper is to prove existence of Beltrami fields with a continuous factor $f$ whose level sets $f^{-1}(k)$ are \textit{nested invariant tori}. The topology of constructed Beltrami fields is consistent with the structure theorem of Arnold. They make a vortex ring that can be supported on a \textit{thick} solid torus.

\vspace{5pt}

%thm1.1
\begin{thm}
For any $u_{\infty}\neq 0$, there exist axisymmetric solutions $u\in C^{1}(\mathbb{R}^{3})$ to (1.4) such that $\nabla \times u$ is compactly supported on a solid torus with some $f\in C(\mathbb{R}^{3})$ such that $f^{-1}(k)$ are nested invariant tori for $0<k< k_0$ with some $k_0$ and empty sets for other $k\neq 0, k_0$. In particular, solutions with regular $f$ exist, e.g. $f\in C^{2+\mu}(\mathbb{R}^{3})$ for $\mu\in (0,1)$. The nested tori are symmetric in the direction of $u_{\infty}$. On each torus, all vortex lines are closed or quasi-periodic. The torus degenerates to a circle for $k=k_0$. 
\end{thm}

\vspace{5pt}

It is worth mentioning that our approach constructing a vortex ring is based on a variational principle for the Grad-Shafranov equation which is available for not only Beltrami flows but also steady Euler flows (1.3) with non-constant $\Pi$ having stream lines and vortex lines not collinear. On the other hand, a vortex can not be supported on a knotted or linked solid torus by this approach due to the axisymmetry.

Vortex rings in Theorem 1.1 are traveling wave solutions to the 3d Euler equations, axisymmetric with swirl. A question on global existence and a finite time blow-up of solutions to the 3d Euler equations may be related with stability and instability of these traveling waves. In the study of shallow water equations, stability of a traveling wave is considered in terms of shape of a wave. Such stability is called \textit{orbital stability}. Orbital stability of a traveling wave of the KdV equation is a well-known property \cite{Ben72}. This property also holds for a traveling wave of the Camassa-Holm equation  \cite{CS00}, \cite{CM01} even though it exhibits a wave breaking. These traveling waves have variational characterization as a minimizer of an energy which deduces their orbital stability by conservation laws.

The Euler equations also admit a variational principle constructing a stationary solution and proving its stability, called a vorticity method, initiated by Arnold \cite[p.88]{AK98}. This method is mainly applied to 2d problems. Benjamin \cite{Ben76} developed it to vortex rings without swirl and suggested their orbital stability, cf. \cite{AC}. A certain class of vortex rings with swirl is constructed by Turkington \cite{Tu89} by a vorticity method. This method involves some quantity \textit{not} conserved by the Euler equations, see Section 2. 

We prove existence of vortex rings by an elementary approach separately from the stability problem by using another variational principle, called a stream function method. This method is a minimax method constructing a vortex ring as a critical point of some functional with given physical constants, e.g. speed and flux of rings. We construct axisymmetric Beltrami fields via the Grad-Shafranov equation as least energy critical points by a minimization on a Nehari manifold.\\

\section{The Grad-Shafranov equation}

\vspace{5pt}

\subsection{Invariant tori}
We derive an equation of vortex rings for the steady Euler equations (1.3). The equation of vortex rings is called Hicks equation, Bragg-Hawthorne equation or Squire-Long equation, e.g. \cite{Fra00}, \cite{Saff}. It is identical to the Grad-Shafranov equation \cite{Grad}, \cite{Shafranov} in magnetohydrodynamics. 

An axisymmetric vector field is written as 

\begin{align*}
u(x)=u^{r}(z,r)e_{r}(\theta)+u^{\theta}(z,r)e_{\theta}(\theta)+u^{z}(z,r)e_{z},
\end{align*}\\
with the cylindrical coordinate $(r,\theta,z)$ and $e_{r}(\theta)={}^{t}(\cos\theta,\sin\theta,0)$, $e_{\theta}(\theta)={}^{t}(-\sin\theta, \cos\theta,0)$, $e_{z}={}^{t}(0,0,1)$, i.e. $x_1=r\cos\theta$, $x_2=r\sin\theta$, $x_3=z$ for $x={}^{t}(x_1,x_2,x_3)$. The each component of $u$ is a function defined in the cross section $\mathbb{R}^{2}_{+}=\{{}^{t}(z,r)\ |\ z\in \mathbb{R},\ r>0 \}$. The divergence-free condition for the axisymmetric vector field is $\partial_z(ru^{z})+\partial_r(ru^{r})=0$. Thus there exists a stream function $\Psi$ such that $ru^{z}=\partial_r\Psi$, $ru^{r}=-\partial_z \Psi$ and $\Psi$ is constant on $\partial\mathbb{R}^{2}_{+}$. We assume that $u\to u_{\infty}$ as $|x|\to\infty$ for $u_{\infty}={}^{t}(0,0,-W)$, $W>0$ and $\Psi(z,0)=-\gamma$ for $\gamma\geq 0$. The circulation of $u$ along the circle $C(z,r)=\{x=re_r(\theta)+ze_z,\ 0\leq \theta\leq 2\pi\}$ is denoted by $2\pi \Gamma(z,r)$. By computation, 

\begin{align*}
\Gamma(z,r)=\frac{1}{2\pi}\int_{C(z,r)}u\cdot \dd l(x)=ru^{\theta}(z,r).
\end{align*}\\
Thus an axisymmetric divergence-free vector field is represented by $\Psi$ and $\Gamma$ as 

\begin{equation*}
\begin{aligned}
u&=\nabla \times \left(\Psi \nabla \theta\right)+\Gamma \nabla \theta. 
\end{aligned}
\tag{2.1}
\end{equation*}\\
This representation is similar as the Clebsch representation, cf. \cite[p.56]{Ben84}, in which the rotation term is replaced with the gradient. The vorticity is written as $\nabla \times u=-L\Psi\nabla \theta + \nabla \Gamma \times \nabla \theta$ with $L=\partial_z^{2}+\partial_r^{2}-r^{-1}\partial_r$.

A stream line of $u$ is a solution to the autonomous system

\begin{align*}
\dot{x}(\tau)=u(x(\tau)).
\end{align*}\\
For 2d incompressible flows, this equation is a (finite-dimensional)  Hamiltonian system with $1$ degree of freedom. Indeed, the velocity is represented by a stream function $\tilde{\Psi}$ as $\nabla_{y}^{\perp}\tilde{\Psi}$ for $\nabla^{\perp}_{y}={}^{t}(\partial_{y_2}, -\partial_{y_1})$ and $y={}^{t}(y_1,y_2)$. Therefore a stream line $y(\tau)$ of 2d flows is a solution to 

\begin{align*}
\dot{y}(\tau)=\nabla^{\perp}_{y}\tilde{\Psi}(y(\tau)).
\end{align*}\\
This is a Hamiltonian system with 1 degree of freedom. Hence a stream line is understood as a level set of the Hamiltonian $\tilde{\Psi}$.

A stream line of axisymmetric Euler flows can be written also by a Hamiltonian system though general 3d flows do not admit this property. We reduce the autonomous system to a Hamiltonian system with 1 degree of freedom, similarly as a Hamiltonian system with 2 degree of freedom with a central field potential, e.g. \cite[p.33]{Arnold}. A stream line $x(\tau)$ of an axisymmetric divergence-free vector field is represented as $x(\tau)=r(\tau)e_{r}(\theta(\tau))+z(\tau)e_{z}$ with the cylindrical coordinate $(r(\tau),\theta(\tau),z(\tau))$. The 2d component $y(\tau)={}^{t}(z(\tau),r(\tau)^{2}/2)$ is a solution to the Hamiltonian system with the stream function $\tilde{\Psi}(y)=\Psi(z,r)$ in the cross section. If $u$ is a solution to the Euler equation (1.3), the circulation $\tilde{\Gamma}(y)=\Gamma(z,r)$ acts as a first integral of the stream line, i.e.  

\begin{align*}
\frac{d}{d\tau}\tilde{\Gamma}(y(\tau))=\dot{z}(\tau)\partial_z\Gamma+\dot{r}(\tau)\partial_r\Gamma=u\cdot \nabla_x \Gamma=0.
\end{align*}\\
Since $y(\tau)$ is a level set of $\tilde{\Psi}$, $\tilde{\Gamma}$ is regarded as a function of $\tilde{\Psi}$, i.e. $\Gamma=\Gamma(\Psi)$. Thus for a given constant $l$, $y(\tau)={}^{t}(z(\tau),r(\tau)^{2}/2)$ is identified with $\tilde{\Psi}^{-1}(l)$ and $\theta(\tau)$ is computed by integrating $\dot{\theta}(\tau)=\Gamma(l)/r(\tau)^{2}$.

If $\tilde{\Psi}^{-1}(l)$ is a closed curve in $\mathbb{R}^{2}_{+}$ for some $l$, the stream line $x(\tau)=r(\tau)e_{r}(\theta(\tau))+z(\tau)e_{z}$ lies on a torus $\Sigma(l)$ in $\mathbb{R}^{3}$. We denote by $\Theta(l)$ increment of $\theta(\tau)$ while $y(\tau)$ goes around the closed curve $\tilde{\Psi}^{-1}(l)$. If $\Theta(l)$ is \textit{commensurable} with $2\pi$, i.e. $\Theta(l)/(2\pi)$ is a rational number, the stream line $x(\tau)$ is closed on the torus. If $\Theta(l)$ is not commensurable with $2\pi$, $x(\tau)$ is quasi-periodic. We call $\Sigma(l)$ an invariant torus \cite{Arnold}.

A vortex line can be computed by the stream line $x(\tau)$. Since the Bernoulli pressure $\tilde{\Pi}(y)=\Pi(z,r)$ is also a first integral of $u$, 

\begin{align*}
\frac{d}{d\tau}\tilde{\Pi}(y(\tau))=\dot{z}(\tau)\partial_z\Pi+\dot{r}(\tau)\partial_r\Pi=u\cdot \nabla_x \Pi=0.
\end{align*}\\
Thus, $\Pi=\Pi(\Psi)$ and the each component of $\omega=\nabla\times u$ is represented as 

\begin{equation*}
\begin{aligned}
\omega^{r}=\dot{\Gamma}(\Psi)u^{r},\quad 
\omega^{\theta}=-r\dot{\Pi}(\Psi)+\dot{\Gamma}(\Psi)u^{\theta},\quad
\omega^{z}=\dot{\Gamma}(\Psi)u^{z}.
\end{aligned}
\tag{2.2}
\end{equation*}\\

\vspace{-10pt}

The Grad-Shafranov equation follows the $\theta$-component of $\omega$ and the boundary conditions for $\Psi$:

\begin{equation*}
\begin{aligned}
-\frac{1}{r^{2}}L\Psi&=-\dot{\Pi}(\Psi)+\frac{1}{r^{2}}\dot{\Gamma}(\Psi)\Gamma(\Psi) \quad \textrm{in}\ \mathbb{R}^2_{+},\\
\Psi&=-\gamma\hspace{95pt} \textrm{on}\ \partial\mathbb{R}^{2}_{+}, \\
\frac{1}{r}\partial_z\Psi&\to 0,\quad \frac{1}{r}\partial_r\Psi\to -W
\hspace{25pt} \textrm{as}\ z^{2}+r^{2}\to \infty.
\end{aligned}
\tag{2.3}
\end{equation*}\\
This equation is a semi-linear elliptic equation for prescribed $\Pi(t)$ and $\Gamma(t)$. It includes axisymmetric flows without swirl $\Gamma\equiv 0$ and (axisymmetric) Beltrami flows $\Pi\equiv \textrm{const}$ with the proportionality factor $f=\dot{\Gamma}(\Psi)$ as particular cases. Strong Beltrami fields can be also described if $\Gamma$ is proportional to $\Psi$. The constants $W>0$ and $\gamma\geq 0$ are speed of a ring and flux measuring distance from the $z$-axis to the ring.

\subsection{The Hicks-Moffatt vortex}
The Grad-Shafranov equation (2.3) admits one explicit solution, called the Hicks-Moffatt vortex, a solution for $\gamma=0$ with 

\begin{align*}
-\dot{\Pi}_{HM}(t)=\lambda_{1}1_{(0,\infty)}(t),\ \lambda_1\in \mathbb{R},\quad \dot{\Gamma}_{HM}(t)=\lambda_2^{1/2}1_{(0,\infty)}(t),\ \lambda_2\geq 0,  
\end{align*}\\
for the indicator function $1_{(0,\infty)}(t)$. The explicit form of $\Psi$ \cite[p.20]{Fra92} is 

\begin{equation*}
\begin{aligned}
\Psi_{HM}(z,r)=
\begin{cases}
&\displaystyle-\frac{3}{2}Wr^{2}\left(B(\kappa)-C(\kappa)\frac{J_{3/2}(\lambda_2^{1/2} \sigma)}{(\lambda_2^{1/2} \sigma)^{3/2}}\right),\quad \sigma< a, \\
&\displaystyle-\frac{1}{2}Wr^{2}\left(1-\frac{a^{3}}{\sigma^{3}}\right),\hspace{83pt} \sigma\geq  a,
\end{cases}
\end{aligned}
\end{equation*}\\
for $\sigma=\sqrt{z^{2}+r^{2}}$, with the constants $a=\kappa/\lambda_2^{1/2}$, 

\begin{align*}
B(\kappa)=\frac{J_{3/2}(\kappa)}{\kappa J_{5/2}(\kappa)},\quad 
C(\kappa)=\frac{\kappa^{1/2}}{J_{5/2}(\kappa)},
\end{align*}\\
and $\kappa\in (0,c_{5/2})$ satisfying 

\begin{align*}
\lambda_1=\frac{3}{2}WB(\kappa)\lambda_2.   
\end{align*}\\
Here, $J_{m}$ denotes the $m$-th order Bessel function of the first kind and $c_{m}$ is the first zero point of $J_{m}$, i.e. $c_{3/2}=4.4934\cdots$, $c_{5/2}=5.7634\cdots$. For given $\lambda_1\in \mathbb{R}$, $\lambda_2>0$ and $W>0$, the constant $\kappa\in (0,c_{5/2})$ uniquely exists since $B(\kappa)$ is decreasing on $(0,c_{5/2})$, i.e. $\kappa=B^{-1}(2\lambda_1/(3W\lambda_2))$. 

The vorticity of the Hicks-Moffatt vortex is supported on a ball $\{x\in \mathbb{R}^{3}\ |\ |x|\leq a \}$ since $\{{}^{t}(z,r)\in \mathbb{R}^{2}_{+}\ |\ \Psi_{HM}>0\}=\{\sigma<a\}$ and $\textrm{spt}\ \omega^{r}, \textrm{spt}\ \omega^{\theta}, \textrm{spt}\ \omega^{z} =\{\sigma\leq a\}$. \\

The Hick-Moffatt vortex is a 3-parameter family of solutions for $\lambda_1\in \mathbb{R}$, $\lambda_2>0$, $W>0$ with the following particular cases:\\

\noindent
\textit{Type I}: $\lambda_2=0$. This is a limiting case. Sending $\lambda_2\to 0$ implies that $\kappa=0$ and

\begin{equation*}
\begin{aligned}
\Psi_{H}(z,r)=
\begin{cases}
&\displaystyle \frac{3}{4}Wr^{2}\left(1-\frac{\sigma^{2}}{a^{2}}\right),\hspace{20pt} \sigma< a, \\
&\displaystyle-\frac{1}{2}Wr^{2}\left(1-\frac{a^{3}}{\sigma^{3}}\right),\hspace{13pt} \sigma\geq a,
\end{cases}
\end{aligned}
\end{equation*}\\
for $\lambda_1>0$ and $W>0$ with $a=(15W/2\lambda_1)^{1/2}$. This is Hill's spherical vortex ring \cite{Hill}.\\

\noindent
\textit{Type II}: $\lambda_1=0$. This implies that $\kappa=c_{3/2}$ and

\begin{equation*}
\begin{aligned}
\Psi_{B}(z,r)=
\begin{cases}
&\displaystyle\frac{3}{2}Wr^{2}C(c_{3/2})\frac{J_{3/2}(\lambda_2^{1/2} \sigma)}{(\lambda_2^{1/2} \sigma)^{3/2}},\quad \sigma< a, \\
&\displaystyle-\frac{1}{2}Wr^{2}\left(1-\frac{a^{3}}{\sigma^{3}}\right),\hspace{44pt} \sigma\geq a,
\end{cases}
\end{aligned}
\end{equation*}\\
for $\lambda_2>0$ and $W>0$ with $a=c_{3/2}/\lambda_2^{1/2}$. The associated velocity $(2.1)$ is a Beltrami field (1.4) with the discontinuous factor $f=\lambda_2^{1/2}1_{(0,\infty)}(\Psi_B)$ \cite{Chandra}. \\

\noindent
\textit{Type III}: $W=0$. This is another limiting case. Sending $W\to 0$ implies that $\kappa=c_{5/2}$ and 

\begin{equation*}
\begin{aligned}
\Psi_{\textrm{com}}(z,r)=
\begin{cases}
&\displaystyle-\frac{\lambda_1}{\lambda_2}r^{2}\left(1-\frac{c_{5/2}^{3/2}}{J_{3/2}(c_{5/2})}\frac{J_{3/2}(\lambda_2^{1/2} \sigma)}{(\lambda_2^{1/2} \sigma)^{3/2}}\right),\quad \sigma< a, \\
&0,\hspace{160pt} \sigma\geq a,
\end{cases}
\end{aligned}
\end{equation*}\\
for $\lambda_1<0$ and $\lambda_2>0$ with $a=c_{5/2}/\lambda_2^{1/2}$. This solution is compactly supported in a half disk $\{{}^{t}(z,r)\in \mathbb{R}^{2}_{+}\ |\ \sigma\leq a\}$ \cite{Prend}.

The Hicks-Moffatt vortex is viewed as a family of solutions for $0<\kappa<c_{5/2}$. The parameter $\lambda_1\in \mathbb{R}$ changes the sign at $\kappa=c_{3/2}$. If $\lambda_1\geq 0$, the Hicks-Moffatt vortex is a unique solution to (2.3) with $\Pi_{HM}(t)$ and $\Gamma_{HM}(t)$ \cite{Fra92}.

\subsection{A free-boundary problem}
A vortex can be supported on a solid torus for $\gamma>0$. The problem (2.3) is a free-boundary problem for $\Psi$ with a priori \textit{unknown} vortex core 

\begin{align*}
\overline{\Omega}=\textrm{spt}\ \omega^{\theta}.
\end{align*}\\
Once the core is found, one can find $\Psi$ by solving two problems:

\begin{align*}
-\frac{1}{r^{2}}L\Psi&=-\dot{\Pi}(\Psi)+\frac{1}{r^{2}}\dot{\Gamma}(\Psi)\Gamma(\Psi) \quad \textrm{in}\ \Omega,\quad \Psi=0\quad \textrm{on}\ \partial\Omega, \\
-\frac{1}{r^{2}}L\Psi&=0 \quad \textrm{in}\ \mathbb{R}^{2}_{+}\backslash \Omega,\quad 
\Psi=-\gamma\quad \textrm{on}\ \partial\mathbb{R}^{2}_{+}, 
\quad \frac{1}{r}\partial_z\Psi\to 0,\hspace{5pt} \frac{1}{r}\partial_r\Psi\to -W
\hspace{5pt} \textrm{as}\ r^{2}+z^{2}\to \infty.
\end{align*}\\
On the other hand, the core is characterized as $\Omega=\{{}^{t}(z,r)\ |\ \Psi>0\  \}$ by a maximum principle if $-L\Psi$ is non-negative. A typical choice of $-\dot{\Pi}(t)$ is a non-negative and non-decreasing function, e.g. 

\begin{align*}
-\dot{\Pi}(t)=\lambda t_{+}^{\alpha},\quad 0\leq \alpha<\infty,\ \lambda>0,
\end{align*}\\
for $t_{+}=\max\{t,0\}$. 

Existence of vortex rings goes back to Fraenkel \cite{Fra70}, \cite{Fra72} who constructed solutions to (2.3) with $\Gamma\equiv 0$ by an implicit function theorem for $\Psi$, called a stream function method. Later, a variational principle is employed to construct solutions for given $W>0$ and $\gamma\geq 0$ with the Lagrange multiplier $\lambda>0$ \cite{FB74}. It is known \cite{AM81}, \cite{Ni80} that solutions can be constructed with given constants $\lambda$, $W$, $\gamma$. 

The another variational principle is a vorticity method of Friedman-Turkington \cite{FT81}, \cite{Friedman82} which provides solutions to (2.3) with $\Gamma\equiv 0$ by maximizing an energy subject to constraints on impulse and mass for $\omega^{\theta}$. The constants $W>0$, $\gamma\geq 0$ are obtained as Lagrange multipliers. The constant $\gamma$ vanishes if impulse is small or mass is large.

Turkington \cite{Tu89} constructed solutions to (2.3) with $\Pi_{HM}(t)$ and $\Gamma_{HM}(t)$ for $\lambda_1 \geq 0$ and $\lambda_2> 0$ with unknown $W>0$, $\gamma\geq 0$ by a vorticity method maximizing

\begin{align*}
\textrm{the kinetic energy}- \frac{1}{2\lambda_2} \int_{\mathbb{R}^{3}}r^{2}\left(\frac{\omega^{\theta}}{r}-\lambda_1\right)_{+}^{2}\dd x,   
\end{align*}\\
subject to constraints on impulse and mass for $\omega^{\theta}$. This result includes the Beltrami flow $\lambda_1=0$. Unlike the case without swirl, the second term is not conserved by the evolution of the Euler equations.

\subsection{Existence of Beltrami fields}

In this paper, we confine ourselves to the problem of Beltrami fields (1.4). The Grad-Shafranov equation (2.3) with $\Pi\equiv \textrm{const}$ is equivalent to (1.4) in the axisymmetric setting with the proportionality factor 

\begin{align*}
f=\dot{\Gamma}(\Psi). 
\end{align*}\\
A level set of $f$ can be identified with that of $\Psi$ if $\dot{\Gamma}(t)$ is invertible. We consider a non-negative and non-decreasing function 

\begin{align*}
\dot{\Gamma}(t)=(\lambda q)^{1/2}t_{+}^{q-1},\quad \lambda > 0,\ 1<  q<\infty.  \tag{2.4}
\end{align*}\\
The function $\dot{\Gamma}(t)$ is continuous at $t=0$ and increasing for $t>0$. The level sets of $f$ are axisymmetric and determined by those of $\Psi$, i.e.  

\begin{equation*}
\begin{aligned}
\textrm{the cross section of}\ f^{-1}(k)=
\begin{cases}
&\left\{{}^{t}(z,r)\in \overline{\mathbb{R}^{2}_{+}}\ \middle|\ \Psi=
\left(\frac{k}{\sqrt{\lambda q}}\right)^{1/(q-1)} \right\}\quad k>0, \\
&\left\{{}^{t}(z,r)\in \overline{\mathbb{R}^{2}_{+}}\ \middle|\ \Psi\leq 0\right\}\hspace{60pt} k=0, \\
&\emptyset \hspace{149pt} k<0.
\end{cases}
\end{aligned}
\end{equation*}\\
If the level set $\Psi^{-1}(l)$ for $l=(k/\sqrt{\lambda q})^{1/(q-1)}$ is a closed curve in $\mathbb{R}^{2}_{+}$, the level set $f^{-1}(k)$ is a torus in $\mathbb{R}^{3}$. By changing the unknown function from $\Psi$ to $\psi$ by setting $\Psi=\psi-Wr^{2}/2-\gamma$, (2.3)-(2.4) with $\Pi\equiv \textrm{const}$ is transformed into the homogeneous problem

\begin{equation*}
\begin{aligned}
-L\psi&=\lambda \left(\psi-\frac{W}{2}r^{2}-\gamma\right)_{+}^{2q-1}\quad \textrm{in}\  \mathbb{R}^2_{+},\\
\psi&=0\hspace{98pt} \textrm{on}\ \partial\mathbb{R}_{+}^{2}, \\
\frac{1}{r}\nabla_{z,r}\psi&\to 0
\hspace{95pt} \textrm{as}\ z^{2}+r^{2}\to \infty.
\end{aligned}
\tag{2.5}
\end{equation*}\\

\vspace{-10pt}

The key result of this paper is existence of solutions to this problem.

\vspace{5pt}

%thm2.1
\begin{thm}
Let $1<q<\infty$ and $0< \lambda, W, \gamma<\infty$. There exists a solution $\psi\in C^{2+\nu}(\overline{\mathbb{R}^{2}_{+}})$ of (2.5) for $\nu\in (0,1)$ such that

\begin{equation*}
\begin{aligned}
&\psi(z,r)=\psi(-z,r),\quad \partial_z\psi(z,r)<0,\  z,r>0,\\
&\Omega=\left\{{}^{t}(z,r)\in \mathbb{R}^{2}_{+}\ \middle|\ \psi(z,r)-\frac{1}{2}Wr^{2}-\gamma>0\right\}.
\end{aligned}
\end{equation*}\\
The vortex core $\Omega$ is bounded, connected and simply-connected with boundary of class $C^{2+\nu}$. The level sets $\{{}^{t}(z,r)\in \mathbb{R}^{2}_{+}\ |\ (\psi-Wr^{2}/2-\gamma)_+=l\}$ are nested closed curves for $0< l< l_0$ with some $l_0>0$, empty sets for other $l\neq 0,l_0$, and points for $l=l_0$.
\end{thm}

\vspace{5pt}

At the limit $q=1$, the function (2.4) is that of the Hicks-Moffatt vortex, i.e. $\dot{\Gamma}_{HM}(t)=\lambda^{1/2}1_{(0,\infty)}(t)$. This function is discontinuous at $t=0$ and constant for $t>0$. A level set of $f=\dot{\Gamma}(\Psi)$ is identified with the vortex core, i.e.

\begin{equation*}
\begin{aligned}
\textrm{the cross section of}\ f^{-1}(k)=
\begin{cases}
&\left\{{}^{t}(z,r)\in \overline{\mathbb{R}^{2}_{+}}\ \middle|\ \Psi>0\right\}\quad k=\lambda^{1/2}, \\
&\left\{{}^{t}(z,r)\in \overline{\mathbb{R}^{2}_{+}}\ \middle|\ \Psi\leq 0\right\}\hspace{11pt} k=0, \\
&\emptyset \hspace{100pt} k\neq \lambda^{1/2}, 0.
\end{cases}
\end{aligned}
\end{equation*}\\
The vortex core of the Hicks-Moffatt vortex (with $\Pi\equiv \textrm{const}$) is a half disk in $\mathbb{R}^{2}_{+}$. Therefore the level set $f^{-1}(k)$ is a ball for $k=\lambda^{1/2}$ in $\mathbb{R}^{3}$. Turkington's solution has a vortex core supported away from the boundary $\partial\mathbb{R}^{2}_{+}$ if $\gamma>0$. Hence $f^{-1}(k)$ is a solid torus for $k=\lambda^{1/2}$ in $\mathbb{R}^{3}$. \\

%section3
\section{A variational principle}

\vspace{5pt}

\subsection{A weighted Sobolev inequality}
In the sequel, we prove Theorem 2.1. Theorem 1.1 is deduced from Theorem 2.1 at the end of this paper. 

Without loss of generality, we may assume that $W=2$ by dividing $\psi$ by $W/2$. We construct a solution to (2.5) by a variational principle. Let $H(\mathbb{R}^{2}_{+}; r^{-1})$ denote the weighted $L^{2}$-Sobolev space of trace zero functions on $\partial\mathbb{R}^{2}_{+}$ with the weight $2\pi^{2}r^{-1}$ \cite[p.243]{OK}. This space is a Hilbert space equipped with the inner product 

\begin{align*}
(\psi,\phi)_{H(\mathbb{R}^{2}_{+}; r^{-1})}=2\pi^{2}\int_{\mathbb{R}^{2}_{+}}\nabla_{z,r} \psi\cdot \nabla_{z,r} \phi \frac{1}{r}\dd z\dd r, 
\end{align*}\\
where $\nabla_{z,r}$ denotes the gradient for the $(z,r)$-variable. The equation (2.5) can be written as a critical point of the functional 

\begin{align*}
I[\psi]=\frac{1}{2}||\psi||_{H(\mathbb{R}^{2}_{+};r^{-1})}^{2}-J[\psi],\quad J[\psi]=\frac{ \pi^{2}\lambda}{q }\int_{\mathbb{R}^{2}_{+}}\left(\psi-r^{2}-\gamma\right)_{+}^{2q}\frac{1}{r}\dd z\dd r.  \tag{3.1} 
\end{align*}\\
This functional is deduced by regarding (2.5) as the 5d problem

\begin{equation*}
\begin{aligned}
-\Delta_{y} \varphi&=\frac{\lambda}{r^{2}}  \left(r^{2}\varphi-r^{2}-\gamma\right)_{+}^{2q-1} \quad \textrm{in}\ \mathbb{R}^{5},
\end{aligned}
\end{equation*}\\ 
by the transform

\begin{align*}
\psi(z,r)\longmapsto
\varphi(y)=\frac{\psi(z,r)}{r^{2}},\quad y={}^{t}(y_1,y')\in \mathbb{R}^{5},\ y_1=z,\ |y'|=r.   
\end{align*}\\

\vspace{-10pt}

The main tool to work with the functional $I$ is the weighted Sobolev inequality \cite[Theorem 2]{Fra92}, 

\begin{align*}
\left(\int_{\mathbb{R}^{2}_{+}}\psi^{p}\frac{1}{r^{2+p/2}}\dd z\dd r\right)^{1/p}
\leq C\left(\int_{\mathbb{R}^{2}_{+}}|\nabla_{z,r} \psi|^{2}\frac{1}{r}\dd z\dd r\right)^{1/2},\quad 2\leq p<\infty.   \tag{3.2}
\end{align*}\\
A general form of this inequality is proved by Koch \cite[Theorem 4.2.2]{Koch99} and de Valeriola-Van Schaftingen \cite[Lemma 3]{Van13}. The inequality (3.2) implies the continuous embedding from $H(\mathbb{R}^{2}_{+};r^{-1})$ to the weighted Lebesgue space 

\begin{align*}
H(\mathbb{R}^{2}_{+};r^{-1})\subset L^{p}(\mathbb{R}^{2}_{+};r^{-2-p/2}),\quad 2\leq p<\infty.    
\end{align*}\\
The functional $I: H(\mathbb{R}^{2}_{+}; r^{-1})\longrightarrow \mathbb{R}$ is bounded since

\begin{align*}
\frac{q}{\pi^{2}\lambda}J[\psi]
\leq \int_{\{\psi>r^{2}\}}\psi^{2q}\frac{1}{r}\dd z\dd r\leq \int_{\mathbb{R}^{2}_{+}}\psi^{p}\frac{1}{r^{2+p/2}}\dd z\dd r\lesssim ||\psi||_{H(\mathbb{R}^{2}_{+};r^{-1})}^{p},\quad p=\frac{2}{3}(4q+1).
\end{align*}\\
It is also continuous by $|s_{+}^{2q}-t_{+}^{2q}|\leq 2q (s_{+}^{2q-1}+t_{+}^{2q-1})((s-t)_{+}+(t-s)_{+})$ and 

\begin{align*}
|J[\psi]-J[\phi]|
\lesssim \left(||\psi||_{H(\mathbb{R}^{2}_{+};r^{-1})}+||\phi||_{H(\mathbb{R}^{2}_{+};r^{-1})} \right)^{p(1-1/(2q))}
||\psi-\phi||_{H(\mathbb{R}^{2}_{+};r^{-1})}^{p/(2q)}.
\end{align*}\\
A similar argument implies that $I$ is Fr\'echet differentiable, i.e. $I\in C^{1}(H(\mathbb{R}^{2}_{+};r^{-1}); \mathbb{R})$, and 

\begin{equation*}
\begin{aligned}
I'[\psi]\phi&=(\psi,\phi)_{H(\mathbb{R}^{2}_{+};r^{-1})}-J'[\psi]\phi,\quad \\ 
J'[\psi]\phi&=2\pi^{2}\lambda \int_{\mathbb{R}^{2}_{+}}(\psi-r^{2}-\gamma)_+^{2q-1}\phi \frac{1}{r}\dd z\dd r,\quad \phi\in H(\mathbb{R}^{2}_{+};r^{-1}). 
\end{aligned}
\tag{3.3}
\end{equation*}\\
A function $\psi \in H(\mathbb{R}^{2}_{+};r^{-1})$ is a critical point $I'[\psi]=0$ if and only if 

\begin{align*}
(\psi,\phi)_{H(\mathbb{R}^{2}_{+};r^{-1})}=2\pi^{2}\lambda\int_{\mathbb{R}^{2}_{+}}(\psi-r^{2}-\gamma)_{+}^{2q-1}\phi\frac{1}{r}\dd z\dd r,\quad \phi\in H(\mathbb{R}^{2}_{+};r^{-1}).     
\end{align*}\\

\vspace{-10pt}

We find a critical point by a minimization on the Nehari manifold 

\begin{align*}
c=\inf_{N(\mathbb{R}^{2}_{+}; r^{-1})}I,\quad N(\mathbb{R}^{2}_{+}; r^{-1})=\left\{\psi\in H(\mathbb{R}^{2}_{+};r^{-1})\ \middle|\ I'[\psi]\psi=0,\ \psi\neq 0 \right\}.     \tag{3.4}
\end{align*}\\
The Nehari manifold includes all critical points of $I$. A minimizer on the Nehari manifold is called a \textit{least energy critical point} or a \textit{grand state} \cite[p.71]{Willem}.

\subsection{Compactness of an embedding}
We construct a grand state in a half space by that in a half disk $D=\{{}^{t}(z,r)\in \mathbb{R}^{2}_{+}\ |\ \sqrt{z^{2}+r^{2}}<R \}$ for $R>0$:  

\begin{equation*}
\begin{aligned}
-L\psi&=\lambda \left(\psi-r^{2}-\gamma\right)_{+}^{2q-1} \quad \textrm{in}\  D,\\ 
 \psi&=0\hspace{85pt} \textrm{on}\ \partial D.
\end{aligned}
\tag{3.5}
\end{equation*}\\ 
A grand state in a half disk is constructed by the same variational principle on the Nehari manifold $N(D;r^{-1})$ and the Hilbert space $H(D;r^{-1})$ defined by the same manner as $\mathbb{R}^{2}_{+}$. By the zero extension, elements of $H(D;r^{-1})$ are regarded as those of $H(\mathbb{R}^{2}_{+};r^{-1})$, i.e. $H(D;r^{-1})\subset H(\mathbb{R}^{2}_{+};r^{-1})$. The functional $I$ is also regarded as that on $H(D;r^{-1})$, i.e. $I\in C^{1}(H(D;r^{-1}); \mathbb{R})$.

The space $H(D;r^{-1})$ is isometrically isomorphic to a homogeneous $L^{2}$-Sobolev space on a ball $B\subset \mathbb{R}^{5}$ centered at the origin with radius $R$ \cite[Lemma 2.2]{AF86}. We denote by $F(B)$ the space of all axisymmetric functions in $H^{1}_{0}(B)$, a homogeneous $L^{2}$-Sobolev space of trace zero functions on $\partial B$ equipped with the inner product $(\varphi,\eta)_{H^{1}_{0}(B)}=\int_{B}\nabla_y \varphi\cdot \nabla_y \eta\dd y$. The transform $\psi(z,r)\longmapsto \varphi(y)=\psi(z,r)/r^{2}$ is a unitary operator from $H(D; r^{-1})$ to $F(B)$, i.e. 

\begin{align*}
(\psi, \phi)_{H(D; r^{-1})}=(\varphi,\eta)_{H^{1}_{0}(B)},\quad \varphi(y)=\frac{\psi(z,r)}{r^{2}},\quad \eta(y)=\frac{\phi(z,r)}{r^{2}}.    \tag{3.6}
\end{align*}\\
Thus $H(D; r^{-1})$ is isometrically isomorphic to $F(B)$, i.e. 

\begin{align*}
H(D; r^{-1})\cong F(B).   
\end{align*}\\

\vspace{-10pt}

The equation for vortex rings of a Beltrami field (3.5) has more \textit{singular} operator and force than those of vortex pairs $-\Delta_{z,r}\psi=\lambda  (\psi-r-\gamma)_{+}^{2q-1}$ \cite{Yang91} and vortex rings without swirl $-L\psi=\lambda r^{2}(\psi-r^{2}-\gamma)_{+}^{2q-1}$ \cite{AM81}, \cite{Ni80}. Associated functionals to these equations admit a grand state by the \textit{compactness} of the embedding $H(D; r^{-1})\subset H^{1}_0(D)\subset \subset L^{p}(D)$ for $p\in [1,\infty)$. Due to the singular force, this compactness is not sufficient to prove existence of a grand state to (3.5).

A key fact to construct a grand state to (3.5) is the following compactness property from $H(D;r^{-1})$ to the weighted Lebesgue space, cf. \cite{OK}. 

\vspace{5pt}

%lem3.1
\begin{lem}[Compact embedding]
\begin{align*}
H(D; r^{-1})\subset \subset L^{p}(D; r^{-1}),\quad 1\leq p<\infty.     \tag{3.7}
\end{align*}
\end{lem}

\vspace{5pt}

\begin{proof}
The continuous embedding follows from (3.2). We prove the compactness for $p=1$. The compactness for $p>1$ follows from the H\"older's inequality. By the Rellich-Kondrachov theorem, 

\begin{align*}
F(B)\subset \subset L^{s}(B),\quad 2<s<\frac{10}{3}.
\end{align*}\\
The conjugate exponent $s'$ satisfies $(-1+3/s)s'<1$. Applying H\"older's inequality implies that $\varphi=\psi/r^{2}$ satisfies 

\begin{align*}
\int_{D}|\psi|\frac{1}{r}\dd z\dd r
\leq \left(\int_{D}|\psi|^{s}\frac{1}{r^{2s-3}}\dd z\dd r\right)^{1/s}\left(\int_{D}\frac{1}{r^{(-1+3/s)s'}}\dd z\dd r\right)^{1/s'} 
\lesssim \left\|\varphi\right\|_{L^{s}(B)}.
\end{align*}\\
Thus, $H(D;r^{-1})\cong F(B)\subset \subset L^{1}(D; r^{-1})$.
\end{proof}

%3.2
\subsection{Regularity of a critical point}
We use a short-hand notation $H=H(D; r^{-1})$. We show that all critical points are classical solutions to (3.5).

\vspace{5pt}

%lem3.2
\begin{lem}
If $\psi\in H$ satisfies $I'[\psi]=0$, then $\varphi=\psi/r^{2}\in C^{2+\nu}(\overline{B})$ for any $\nu\in (0,1)$, and  

\begin{equation*}
\begin{aligned}
-\Delta_y \varphi&=\frac{\lambda}{r^{2}}(r^{2}\varphi-r^{2}-\gamma)^{2q-1}_{+}\quad \textrm{in}\ B,\\
\varphi&=0\hspace{96pt} \textrm{on}\ \partial B. 
\end{aligned}
\tag{3.8}
\end{equation*}\\
In particular, $\psi\in C^{2+\nu}(\overline{D})$ is a classical solution to (3.5).
\end{lem}

\vspace{5pt}

\begin{proof}
The force term of the 5d problem for $\varphi$ belongs to $L^{3}(B)$ by (3.2) and 

\begin{align*}
\int_{B}\left|\frac{\lambda}{r^{2}}(r^{2}\varphi-r^{2}-\gamma)_{+}^{2q-1}  \right|^{3}\dd y 
&=2\pi^{2}\lambda^{3}\int_{D}\frac{1}{r^{3}}(\psi-r^{2}-\gamma)_{+}^{3(2q-1)}\dd z\dd r  \\
&\leq 2\pi^{2}\lambda^{3}\int_{D}|\psi|^{p}\frac{1}{r^{2+p/2}}\dd z\dd r\lesssim ||\psi||_{H}^{p},\quad p=8q-\frac{14}{3}.
\end{align*}\\
Thus an axisymmetric solution of the Poisson equation $-\Delta_y\tilde{\varphi}=\lambda r^{-2}(r^{2}\varphi-r^{2}-\gamma)_{+}^{2q-1}$ in $B$ and $\tilde{\varphi}=0$ on $\partial B$ belongs to $W^{2,3}(B)$. Here, $W^{k,s}(B)$ denotes the Sobolev space of order $k$ with exponent $s$. Since $I'[\psi]=0$ and $\tilde{\psi}=r^{2}\tilde{\varphi}$ satisfies  

\begin{align*}
(\tilde{\psi},\phi)_{H}=2\pi^{2}\lambda\int_{D}(\psi-r^{2}-\gamma)_{+}^{2q-1}\phi\frac{1}{r}\dd z\dd r,\quad \phi\in H,
\end{align*}\\
the function $\tilde{\psi}$ agrees with $\psi$. Thus $\varphi=\tilde{\varphi}\in W^{2,3}(B)$ is a solution to (3.8). By the Sobolev inequality, $\varphi\in L^{\infty}(B)$ and $-\Delta_y \varphi\in L^{\infty}(B)$. Thus $\varphi\in W^{2,s}(B)$, $1\leq s<\infty$, by elliptic regularity. Since $\gamma>0$,  

\begin{align*}
-\Delta_y\varphi=\lambda r^{4(q-1)}\left(\varphi-1-\frac{\gamma}{r^{2}}\right)_{+}^{2q-1}\in W^{1,\infty}(B).
\end{align*}\\
Thus $\varphi\in W^{3,s}(B)$ and $\varphi\in C^{2+\nu}(\overline{B})$ for any $\nu\in (0,1)$.\\
\end{proof}

%section3
\section{Grand states in a half disk}

\vspace{5pt}

\subsection{A deformation theorem}
In this section, we prove existence of a grand state and its properties in a half disk. We consider the minimization

\begin{align*}
c=\inf_{ N} I ,\quad
N=\left\{\psi\in H\ \middle|\ I'[\psi]\psi=0,\ \psi\neq 0   \right\}.    \tag{4.1}
\end{align*}\\
By the compact embedding $H=H(D;r^{-1})\subset\subset L^{p}(D;r^{-1})$ for $p\in [1,\infty)$, the functional $I$ satisfies the Palais-Smale condition and admits a deformation theorem.

\vspace{5pt}

%prop3.1
\begin{prop}
\begin{align*}
\frac{1}{2}\left(1-\frac{1}{q}\right)||\psi||_{H}^{2}\leq I[\psi]-\frac{1}{2q}I'[\psi]\psi,    \quad \psi\in H.       \tag{4.2}
\end{align*}
\end{prop}

\vspace{5pt}

\begin{proof}
This follows from (3.1) and (3.3). 
\end{proof}

\vspace{5pt}

%prop3.2
\begin{prop}[Palais-Smale condition]
Any sequence $\{\psi_n\}\subset H$ satisfying 

\begin{align*}
\sup_{n}I[\psi_n]<\infty,\quad ||I'[\psi_n]||_{H^{*}}\to 0,  
\end{align*}\\
has a convergent subsequence in $H$. 
\end{prop}

\vspace{5pt}

\begin{proof}
A sequence satisfying the above condition is bounded in $H$ by (4.2) and has a convergent subsequence weakly in $H=H(D; r^{-1})$ and strongly in $L^{p}(D; r^{-1})$ for $1\leq p<\infty$ by (3.7). Thus by choosing a subsequence (still denoted by $\{\psi_n\}$), there exists some $\psi$ such that 

\begin{equation*}
\begin{aligned}
\psi_{n}&\rightharpoonup \psi \quad \textrm{in}\ H(D; r^{-1}),\\
\psi_n&\to \psi\quad \textrm{in}\ L^{p}(D; r^{-1}).
\end{aligned}
\end{equation*}\\
This implies that $J'[\psi_n]\psi_n$ and $J'[\psi_n]\psi$ converge to $J'[\psi]\psi$. By (3.3), 

\begin{align*}
I'[\psi_n]\psi&=(\psi_n,\psi)_{H}-J'[\psi_n]\psi,\\
I'[\psi_n]\psi_n&=(\psi_n,\psi_n)_{H}-J'[\psi_n]\psi_n.
\end{align*}\\
Since the left-hand sides vanish as $n\to\infty$, $\lim_{n\to\infty}||\psi_n||_{H}=||\psi||_{H}$. Hence $\psi_n\to\psi$ in $H$.
\end{proof}

\vspace{5pt}

We set a filtration and a set of critical points with a critical value $c\in \mathbb{R}$ by 

\begin{equation*}
\begin{aligned}
A_c&=\{\psi\in H\ |\ I[\psi]\leq c\ \}, \\
K_c&=\{\psi\in H\ |\ I[\psi]=c,\ I'[\psi]=0\ \}.
\end{aligned}
\end{equation*}

\vspace{5pt}

%lem3.3
\begin{lem}[Deformation Theorem]
There exists $\varepsilon_0>0$ such that for $c\in \mathbb{R}$ and a neighborhood ${\mathcal{U}}$ of $K_c$, there exists $\varepsilon_1\in (0,\varepsilon_0)$ and a homeomorphism $\iota: H\to H$ such that 

\begin{equation*}
\begin{aligned}
&\iota(\psi)=\psi,\quad \psi\notin I^{-1}[c-\varepsilon_0,c+\varepsilon_0],   \\
&I[\iota(\psi)]\leq I[\psi], \quad \psi\in H,  \\
&\iota(A_{c+\varepsilon_1}\backslash {\mathcal{U}})\subset A_{c-\varepsilon_1}.   
\end{aligned}
\end{equation*}
\end{lem}

\vspace{5pt}

\begin{proof}
See \cite[p.82, Theorem A.4]{Rab}.
\end{proof}

%3.2
\subsection{Characterization of a critical value}

We characterize a critical value (4.1) as a minimax value by using convexity of $I[t\psi]$ for $t>0$. Then the deformation theorem implies that a grand state is a critical point, cf. \cite[p.74, Theorem 4.3]{Willem}.

\vspace{5pt}

%prop4.4
\begin{prop}
For $\psi\in H\backslash \{0\}$, set 

\begin{align*}
g(t)=I[t\psi],\quad t\geq 0.   
\end{align*}\\
There exists some $t(\psi)>0$ such that $\dot{g}(t(\psi))=0$, i.e. $t(\psi)\psi\in N$. The function $g(t)$ is increasing for $t<t(\psi)$ and decreasing for  $t>t(\psi)$. Moreover, 

\begin{align*}
g(t(\psi))=\frac{\pi^{2}\lambda}{q}\int_{D}(t(\psi)\psi-r^{2}-\gamma )_{+}^{2q-1}\left(r^{2}+\gamma+(q-1)t(\psi)\psi \right)\frac{1}{r}\dd z\dd r, \tag{4.3}
\end{align*}\\
and $t(\cdot): H\backslash \{0\}\to (0,\infty)$ is continuous.
\end{prop}

\vspace{5pt}

\begin{proof}
By

\begin{align*}
&\frac{\dot{g}(t)}{t}=||\psi||_{H}^{2}-\frac{2\pi^{2}\lambda}{t}\int_{D}(t\psi-r^{2}-\gamma)_{+}^{2q-1}\psi \frac{1}{r}\dd z\dd r, \\
&\frac{\dd}{\dd t}\left(\frac{\dot{g}(t)}{t}\right)
=-\frac{2\pi^{2}\lambda}{t^{2}}\int_{D}(t\psi-r^{2}-\gamma)_{+}^{2q-2}\left( 2(q-1)t\psi+r^{2}+\gamma  \right)\psi \frac{1}{r}\dd z\dd r<0,
\end{align*}\\
$\lim_{t\to 0}\dot{g}(t)/t=||\psi||_{H}^{2}>0$ and $\dot{g}(t)/t$ is decreasing. Hence there exists a unique $t(\psi)>0$ such that $\dot{g}(t(\psi))=0$ and $t(\psi)\psi\in N$. The identity (4.3) follows from 

\begin{align*}
g(t(\psi))=I[t(\psi)\psi]=\frac{1}{2}||t(\psi)\psi||_{H}^{2}-J[t(\psi)\psi]=\frac{1}{2}J'[t(\psi)\psi](t(\psi)\psi)-J[t(\psi)\psi]. 
\end{align*}\\
To prove continuity of $t(\cdot)$, we take a sequence $\{\psi_n\}\subset H\backslash \{0\}$ such that $\psi_n\to \psi$ in $H\backslash \{0\}$. Since $g_n(t)=I[t\psi_n]$ satisfies $\dot{g}_n(t(\psi_n))=0$, 

\begin{align*}
0=\frac{\dot{g}_{n}(t(\psi_n))}{t(\psi_n)}=||\psi_n||_{H}^{2}-2\pi^{2}\lambda t(\psi_n)^{2(q-1)}\int_{D}\left(\psi_n-\frac{r^{2}+\gamma}{t(\psi_n)}\right)_{+}^{2q-1}\psi_n \frac{1}{r}\dd z\dd r.
\end{align*}\\
Since $\lim_{n\to\infty}||\psi_n||_{H}=||\psi||_{H}\neq 0$, the sequence $\{t(\psi_n)\}$ is bounded. 

Suppose that $\{t(\psi_n)\}$ does not converge to $t(\psi)$. Then, there exists a subsequence (still denoted by $\{t(\psi_n)\}$) such that $t(\psi_n)\to t_0$ for some $t_0\geq 0$. Since $\psi\neq 0$, $t_0>0$. Sending $n\to\infty$ to the above equality implies that $\dot{g}(t_0)=0$ for $g(t)=I[t\psi]$. Thus $t_0=t(\psi)$. This is a contradiction and we conclude that $t(\psi_n)\to t(\psi)$.
\end{proof}

\vspace{5pt}

%prop4.5
\begin{prop}
\begin{align*}
c=\inf_{\psi\in N}I[\psi]=\inf_{\psi\in H\backslash \{0\}}\sup_{t\geq 0}I[t\psi]=\inf_{p\in \Lambda}\sup_{0\leq t\leq 1}I[p(t)],   \tag{4.4}
\end{align*}\\
for $\Lambda=\{p\in C([0,1]; H)\ |\ p(0)=0,\ I[p(1)]<0\ \}$.
\end{prop}

\vspace{5pt}

\begin{proof}
Since $\sup_{t\geq 0}I[t\psi]\leq I[\psi]$ for $\psi\in N$ by Proposition 4.4, 

\begin{align*}
\inf_{\psi\in H\backslash \{0\}}\sup_{t\geq 0}I[t\psi]\leq \inf_{\psi\in N}I[\psi]=c.
\end{align*}\\
Since 

\begin{align*}
\inf_{p\in \Lambda}\sup_{0\leq t\leq 1}I[p(t)]\leq \inf_{\psi\in H\backslash \{0\}}\sup_{t\geq 0}I[t\psi], 
\end{align*}\\
it suffices to show that the left-hand side is larger than $c$. We set $h(t)=I'[p(t)]p(t)$ for $p\in \Lambda$. By (3.3) and (3.7),

\begin{align*}
I'[\psi]\psi
=||\psi||_{H}^{2}-2\pi^{2}\lambda\int_{D}(\psi-r^{2}-\gamma)_{+}^{2q-1}\psi\frac{1}{r}\dd z\dd r
\geq ||\psi||_{H}^{2}-C||\psi||_{H}^{2q},   \quad \psi\in H,
\end{align*}\\
for some $C>0$. This implies that $\lim_{t\to 0}h(t)/||p(t)||_{H}^{2}=1$. Hence $h(t)$ is positive near $t=0$. Since $h(1)<0$ by (4.2), by the intermediate value theorem, $h(t_1)=0$ for some $t_1\in (0,1)$. Thus $p(t_1)\in N$ and

\begin{align*}
c=\inf_{\psi \in N }I[\psi]\leq I[p(t_1)]\leq \sup_{0\leq t\leq 1}I[p(t)].
\end{align*}\\
Since $p\in \Lambda$ is arbitrary, (4.4) follows.
\end{proof}

\vspace{5pt}

%lem4.6
\begin{lem}
If $\psi\in N$ satisfies $c=I[\psi]$, then $I'[\psi]=0$.
\end{lem}

\vspace{5pt}

\begin{proof}
Suppose on the contrary that there exists a grand state $\psi\in N$ such that $I'[\psi]\neq 0$. By the continuity of $I': H\to H^{*}$, there exists $\delta>0$ such that $I'[\phi]\neq 0$ for all $\phi\in \mathcal{B}$, where $\mathcal{B}=\{\phi\in H\ |\ ||\psi-\phi||_{H}\leq \delta \}$. Thus the set of critical points with the critical value $c$ is not included in $\mathcal{B}$, i.e. $K_{c}\cap \mathcal{B}=\emptyset$. Thus $\mathcal{U}=\mathcal{B}^{c}$ is a neighborhood of $K_{c}$. We apply Lemma 4.3 and take $\varepsilon_1>0$ and a homeomorphism $\iota: H\to H$ such that

\begin{align*}
\iota(A_{c+\varepsilon_1}\cap \mathcal{B} )\subset A_{c-\varepsilon_1}.
\end{align*}\\
We take $p\in \Lambda$ such that $p(t)\in A_{c-\varepsilon_1}\cup( A_{c+\varepsilon_1}\cap \mathcal{B})$ for all $0\leq t\leq 1$ and set $\hat{p}=\iota(p)\in \Lambda$. Then, $\hat{p}(t)\in A_{c-\varepsilon_1}$ for all $0\leq t\leq 1$. Thus by (4.4),  

\begin{align*}
c\leq \sup_{0\leq t\leq 1}I[\tilde{p}(t)]\leq c-\varepsilon_1.
\end{align*}\\
This is a contradiction. We conclude that $I'[\psi]=0$.
\end{proof}

%3.3
\subsection{Existence of a grand state}
We prove existence of a symmetric grand state for the $z$-variable. We say that $\psi^{*}$ is the Steiner symmetrization of $\psi$ if $\psi^{*}(z,r)=\psi^{*}(-z,r)$, ${}^{t}(z,r)\in D$, $\psi^{*}$ is non-increasing for $|z|$, and $\psi^{*}$ is equi-mesurable, i.e.

\begin{align*}
|\{z\in \mathbb{R}\ |\ \psi(z,r)\geq t,\ {}^{t}(z,r)\in D \}|=|\{z\in \mathbb{R}\ |\ \psi^{*}(z,r)\geq t,\ {}^{t}(z,r)\in D \}|,\quad t\geq 0,\ r\geq 0.
\end{align*}\\
The Steiner symmetrization exists for any $\psi\in H$ and does not increase the Dirichlet energy \cite[Appendix I]{FB74}, i.e. $||\psi^{*}||_{H}\leq ||\psi||_{H}$. We show that a grand state can be replaced with the Steiner symmetrization.

\vspace{5pt}

%4.7
\begin{prop}
If $\psi\in N$ satisfies $c=I[\psi]$, then $\psi^{*}\in N$ and $c=I[\psi^{*}]$.
\end{prop}

\vspace{5pt}

\begin{proof}
Since $\psi^{*}$ is equi-mesurable, 

\begin{align*}
I'[\psi^{*}]\psi^{*}
&=||\psi^{*}||_{H}^{2}-2\pi^{2}\lambda \int_{D}(\psi^{*}-r^{2}-\gamma)_{+}^{2q-1}\psi^{*}\frac{1}{r}\dd z\dd r \\
&\leq ||\psi||_{H}^{2}-2\pi^{2}\lambda \int_{D}(\psi-r^{2}-\gamma)_{+}^{2q-1}\psi\frac{1}{r}\dd z\dd r=I'[\psi]\psi=0.
\end{align*}\\
Thus $g(t)=I[t\psi^{*}]$ satisfies $\dot{g}(1)\leq 0$. Since $g(t)$ is decreasing for $t>t(\psi^{*})$ by Proposition 4.4, we have $t(\psi^{*})\leq 1$. By (4.3), $g(t(\psi^{*}))\leq g(1)$ and 

\begin{align*}
c\leq I[t(\psi^{*})\psi^{*}]=g(t(\psi^{*}))\leq g(1)=I[\psi^{*}]\leq I[\psi]=c. 
\end{align*}\\
Thus $t(\psi^{*})=1$, $\psi^{*}\in N$ and $c=I[\psi^{*}]$.
\end{proof}

\vspace{5pt}

%prop4.8
\begin{prop}
\begin{align*}
\inf_{\psi\in N}||\psi||_{H}>0.    \tag{4.5}
\end{align*}
\end{prop}

\vspace{5pt}

\begin{proof}
For $\varepsilon>0$, we take $\delta>0$ such that $|s|^{2(q-1)}\leq \varepsilon$ for $|s|\leq \delta$. Applying (3.7) to $\psi\in N$ yields 

\begin{align*}
||\psi||_{H}^{2}=J'[\psi]\psi
=\int_{\{\psi<\delta\}}+\int_{\{\psi\geq \delta\}} 
\lesssim \varepsilon\int_{D}\psi^{2}\frac{1}{r}\dd z\dd r +\int_{D}\psi^{2q}\frac{1}{r}\dd z\dd r 
\lesssim \varepsilon ||\psi||_{H}^{2}+||\psi||_{H}^{2q}.
 \end{align*}\\
Thus the desired result follows.
\end{proof}

\vspace{5pt}

%lem4.9
\begin{lem}[Existence of a grand state]
There exists $\psi\in N$ such that $c=I[\psi]$, $I'[\psi]=0$ and $\psi=\psi^{*}$.
\end{lem}

\vspace{5pt}

\begin{proof}
We take a minimizing sequence $\{\psi_n\}\subset N$ of (4.1). By (4.2) and (3.7), there exists a subsequence (still denoted by $\{\psi_n\}$) and some $\psi$ such that 

\begin{equation*}
\begin{aligned}
\psi_{n}&\rightharpoonup \psi \quad \textrm{in}\ H(D; r^{-1}),\\
\psi_n&\to \psi\quad \textrm{in}\ L^{p}(D; r^{-1}),\quad 1\leq p<\infty.
\end{aligned}
\end{equation*}\\
The limit $\psi$ is non-trivial by (4.5). This convergence implies $J'[\psi_n]\psi_n\to J'[\psi]\psi$ and 

\begin{align*}
||\psi||_{H}^{2}-J'[\psi]\psi\leq\liminf_{n\to\infty}\left(||\psi_n||_{H}^{2}-J'[\psi_n]\psi_n\right)=\lim_{n\to\infty}I'[\psi_n]\psi_n=0.  
\end{align*}\\
Hence $g(t)=I[t\psi]$ satisfies $\dot{g}(1)\leq 0$. Since $g(t)$ is decreasing for $t>t(\psi)$ for some $t(\psi)>0$ by Proposition 4.4, we have $t(\psi)\leq 1$. By (4.3), $g(t(\psi))\leq g(1)$ and 

\begin{align*}
\inf_{N}I\leq g(t(\psi))\leq g(1)\leq\lim_{n\to\infty}I[\psi_n]= \inf_{N}I.
\end{align*}\\
Hence, $t(\psi)=1$, $\psi\in N$ and $c=I[\psi]$. Thus $\psi$ is a grand state. By Proposition 4.7, we replace it by the Steiner symmetrization.  
\end{proof}

%3.4
\subsection{Shape of vortex cores}

We set $\Omega=\{{}^{t}(z,r)\in D\ |\ \psi-r^{2}-\gamma>0\}$ for the grand state $\psi$. By (2.2) and (2.4), $\Omega$ is the vortex core, i.e. $\overline{\Omega}=\textrm{spt}\ \omega^{\theta}$. The vortex core $\Omega$ is symmetric in the $z$-direction since $\psi$ is. We show that $\psi$ is decreasing for $|z|$ and $\Omega$ consists of simply-connected components with regular boundaries, cf. \cite[Theorem 3D]{FB74}.

\vspace{5pt}

%prop
\begin{prop}
The grand state in Lemma 4.9 satisfies $\psi\in C^{2+\nu}(\overline{D})$ for $\nu\in (0,1)$ and 

\begin{align*}
\frac{\partial \psi}{\partial z}(z,r)<0,\quad z>0,\ {}^{t}(z,r)\in D.   \tag{4.6}
\end{align*}\\
The level sets $\{\psi-r^{2}-\gamma=l\}$ are nested closed curves of class $C^{2+\nu}$ for $0\leq l<l_0$ with $l_0=||(\psi-r^{2}-\gamma)_{+}||_{\infty}$ and points for $l=l_0$.
\end{prop}

\vspace{5pt}

\begin{proof}
The regularity of the grand state follows from Lemma 3.2. Since $\psi$ is Steiner symmetric,  

\begin{align*}
\frac{\partial \psi}{\partial z}(z,r)\leq 0,\quad z>0,\ {}^{t}(z,r)\in D.    
\end{align*}\\
We set $\varphi(y)=\psi(z,r)/r^{2}$ for $y={}^{t}(y_1,y')$, $y_1=z$, $|y'|=r$. By Lemma 3.2, $\varphi\in C^{2+\nu}(\overline{B})$ for $\nu\in (0,1)$ and $\varphi$ is a solution to the 5d Dirichlet problem (3.8) in $B$. For arbitrary $\tau\in (0,R)$, we set a hyperplane and a cap by 

\begin{align*}
T_{\tau}=\{{}^{t}(y_1,y')\in B\ |\ y_1=\tau  \},\quad \Sigma_{\tau}=\{{}^{t}(y_1,y')\in B\ |\ y_1>\tau  \}.
\end{align*}\\
We denote a reflection point of $y={}^{t}(y_1,y')\in \Sigma_{\tau}$ with respect to $T_{\tau}$ by $y_{\tau}={}^{t}(2\tau-y_1, y')$ and set $\tilde{\varphi}(y)=\varphi(y)-\varphi(y_{\tau})$ for $y\in \Sigma_{\tau}$. Since $\varphi(y_{\kappa})$ solves (3.8) for $y\in \Sigma_{\kappa}$ and $\varphi(y)$ is non-increasing for $y_1>0$, $\tilde{\varphi}(y)\leq 0$ in $\Sigma_{\tau}$ and 

\begin{align*}
-\Delta_y \tilde{\varphi}(y)&=\frac{\lambda}{r^{2}}(r^{2}\varphi(y)-r^{2}-\gamma)_{+}^{2q-1}-\frac{\lambda}{r^{2}}(r^{2}\varphi(y_{\kappa})-r^{2}-\gamma)_{+}^{2q-1}\leq 0,\quad y\in \Sigma_{\kappa},\\
\tilde{\varphi}(y)&=0,\quad y\in  T_{\kappa}.
\end{align*}\\
Hence by Hopf's lemma \cite[p.65, Theorem 7]{PW}, 

\begin{align*}
\frac{\partial \tilde{\varphi}}{\partial y_1}(y)<0 ,\quad  y\in T_{\tau}.
\end{align*}\\
Thus (4.6) holds. By the implicit function theorem, a level set of $\psi-r^{2}-\gamma=l$ for $0\leq l< l_0$ is written as a graph of a $C^{2+\nu}$-function. For $l=l_0$, the level set is points lying on the $r$-axis.
\end{proof}

\vspace{5pt}

The connectedness of $\Omega$ follows from the least energy property of a grand state, cf. \cite[Theorem 4]{AM81}. Our proof is based on that of de Valeriola-Van Schaftingen \cite[Lemma 11]{Van13} using energy identities  \cite[Lemma 5A]{FB74}.

\vspace{5pt}

%prop4.11
\begin{prop}
The identities 
\begin{align*}
\int_{\Omega}|\nabla_{z,r} \Psi|^{2}&=\lambda \int_{\Omega}\Psi^{2q},   \tag{4.7}\\
\int_{\Omega}|\nabla_{z,r} \Psi|^{2}&= \int_{\Omega}|\nabla_{z,r} \psi|^{2}-\int_{\Omega}|\nabla_{z,r} (r^{2}+\gamma)|^{2},   \tag{4.8}
\end{align*}\\
hold for $\Psi=\psi-r^{2}-\gamma$ and the grand state $\psi$ in Lemma 4.9, where the measure $2\pi^{2}r^{-1}\dd z\dd r$ is suppressed.
\end{prop}

\vspace{5pt}

\begin{proof}
Since $\Psi$ satisfies $-L \Psi=\lambda \Psi^{2q-1}$ in $\Omega$ and $\Psi=0$ on $\partial\Omega$, $\Phi=\Psi/r^{2}$ satisfies the 5d problem 

\begin{align*}
-\Delta_y \Phi=\lambda r^{4(q-1)}\Phi^{2q-1}\quad \textrm{in}\ U,\qquad 
\Phi=0\quad \textrm{on}\ \partial U,
\end{align*}\\
for $U=\{y={}^{t}(y_1,y')\in B\ |\  y_1=z,\ |y'|=r,\ {}^{t}(z,r)\in \Omega  \}$. By multiplying $\Phi$ by the equation and integration by parts,

\begin{align*}
\int_{U}|\nabla_{y} \Phi|^{2}\dd y=\lambda\int_{U}r^{4(q-1)}\Phi^{2q}\dd y.
\end{align*}\\
Since $H(D;r^{-1})\cong F(B)$ by the transform $\Psi\longmapsto \Phi=\Psi/r^{2}$, (4.7) follows.

The identity (4.8) follows from $|\nabla_{z,r} \Psi|^{2}
=|\nabla_{z,r} \psi|^{2}-|\nabla_{z,r} (r^{2}+\gamma)|^{2}
-2\nabla_{z,r} \Psi\cdot \nabla_{z,r} (r^{2}+\gamma)$, and 

\begin{align*}
\int_{\Omega}\nabla_{z,r} \Psi\cdot \nabla_{z,r} (r^{2}+\gamma)
=\int_{U}\nabla_y \Phi\cdot \nabla_y \left(1+\frac{\gamma}{r^{2}}\right)\dd y=0.
\end{align*}
\end{proof}

\vspace{5pt}

%lem3.12
\begin{lem}
The vortex core $\Omega$ of the grand state in Lemma 4.9 is connected.
\end{lem}

\vspace{5pt}

\begin{proof}
Let $\Omega_0 \subset \Omega$ be a connected component of $\Omega=\Omega_0\cup \Omega_1$. We set 

\begin{align*}
\Psi_0=
\begin{cases}
\ \Psi & \textrm{in}\ \Omega_0, \\
\ 0 & \textrm{in}\ D\backslash \overline{\Omega_0},
\end{cases}
\end{align*}\\
and set $\psi_0=\Psi_0+\chi$ with $\chi=\min\{\psi,r^{2}+\gamma  \}$. Then,

\begin{equation*}
\begin{aligned}
\psi_0=
\begin{cases}
\ \psi &  \textrm{in}\ D\backslash \overline{\Omega_1}, \\
\ r^{2}+\gamma &  \textrm{in}\ \Omega_1.
\end{cases}
\end{aligned}
\end{equation*}\\
By $I'[\psi]\psi=0$, 

\begin{align*}
I'[\psi_0]\psi_0
&=||\psi_0||_{H}^{2}-\lambda\int_{D}(\psi_0-r^{2}-\gamma)_{+}^{2q-1}\psi_0 \\
&=-\int_{\Omega}|\nabla_{z,r} \psi|^{2}+\int_{\Omega}|\nabla_{z,r} \chi|^{2}
+\int_{\Omega_0}|\nabla_{z,r} \psi|^{2}-\int_{\Omega_0}|\nabla_{z,r} \chi|^{2}
+\lambda \int_{\Omega_1}\Psi_{+}^{2q-1}\psi.
\end{align*}\\
Since the identities (4.7) and (4.8) hold also on $\Omega_0$, 

\begin{align*}
I'[\psi_0]\psi_0=\lambda\int_{\Omega_1}\Psi_+^{2q-1}\chi\geq 0.
\end{align*}\\
Thus $I'[t\Psi_0+\chi](t\Psi_0+\chi)$ is non-negative at $t=1$. By (3.3), $I'[t\Psi_0+\chi](t\Psi_0+\chi)$ is negative for large $t>1$. By the intermediate value theorem, there exists $t_0\geq 1$ such that $t_0\Psi_0+\chi\in N$. By (4.8), 

\begin{align*}
\int_{D}|\nabla_{z,r} \chi|^{2}=\int_{D\backslash \Omega}|\nabla_{z,r} \psi|^{2}+\int_{\Omega}|\nabla_{z,r} (r^{2}+\gamma)|^{2}
=\int_{D}|\nabla_{z,r} \psi|^{2}-\int_{\Omega}|\nabla_{z,r} \Psi|^{2}. 
\end{align*}\\
Applying (4.7) yields

\begin{align*}
I[\psi]\leq 
I[t_0\Psi_0+\chi]
&=\frac{1}{2}||t_0\Psi_0+\chi ||_{H}^{2}
-\frac{\lambda}{2q} \int_{D}(t_0\Psi_0+\chi-r^{2}-\gamma)_{+}^{2q}  \\
&=I[\psi]+\frac{t_0^{2}}{2}\left(1-\frac{t_0^{2(q-1)}}{q} \right)\int_{\Omega_0}|\nabla_{z,r} \Psi|^{2}-\frac{1}{2}\left(1-\frac{1}{q}\right)\int_{\Omega}|\nabla_{z,r} \Psi|^{2} \\
&\leq I[\psi]-\frac{1}{2}\left(1-\frac{1}{q}\right)\int_{\Omega_1}|\nabla_{z,r} \Psi|^{2} \leq I[\psi]. 
\end{align*}\\
We conclude that $\Omega_0=\Omega$.\\
\end{proof}

%subsection3.5
\section{Grand states in a half plane}

\vspace{5pt}

%5.1
\subsection{A pointwise estimate of a grand state}

We construct a grand state in $\mathbb{R}^{2}_{+}$ by sending $R\to\infty$ to that in a half disk $D=D(R)$. To take a limit, we estimate a grand state in  $L^{\infty}$ uniformly for $R$ by using the Green function of the Dirichlet problem 

\begin{equation*}
\begin{aligned}
-L\psi&=r^{2}\zeta\quad \textrm{in}\ \mathbb{R}^{2}_{+}, \\
\psi&=0\hspace{20pt} \textrm{on}\ \partial\mathbb{R}^{2}_{+}.
\end{aligned}
\end{equation*}\\
The Grad-Shafranov equation (2.5) is the problem for $\zeta=\lambda r^{-2}(\psi-r^{2}-\gamma)_{+}^{2q-1}$. The operator $-r^{-2}L$ is viewed as the 5d Laplace operator by the transform $\psi\longmapsto \varphi=\psi/r^{2}$, i.e. $-\Delta_{y}\varphi=\zeta$ in $\mathbb{R}^{5}$. Solutions to this problem are  represented by 

\begin{align*}
\psi(z,r)=\int_{\mathbb{R}^{2}_{+}}G(z,r,z',r')\zeta(z',r')r'\dd z'\dd r',  \tag{5.1}
\end{align*}\\
with the Green function  

\begin{align*}
G(z,r,z',r')= \frac{rr'}{2\pi}\int_{0}^{\pi}\frac{\cos\theta \dd \theta}{\sqrt{|z-z'|^{2}+|r-r'|^{2}-2rr'\cos\theta  }}.
\end{align*}\\
The Green function is written by the complete elliptic integrals of the first and second kind. By their asymptotic expansions with $\xi^{2}= (|z-z'|^{2}+|r-r'|^{2}) /{4rr'}$, e.g. \cite[p.482]{Friedman82},  

\begin{align*}
\frac{G(z,r,z',r')}{\sqrt{rr'}}\lesssim
\begin{cases}
&|\log\xi| \quad \xi\leq 1, \\
& \xi^{-3}\hspace{24pt} \xi\geq 1.
\end{cases}
\end{align*}\\
Thus the Green function satisfies the pointwise estimate

\begin{align*}
G(z,r,z',r')\leq C\frac{(rr')^{1/2+\tau} }{(|z-z'|^{2}+|r-r'|^{2})^{\tau} },\quad 0<\tau\leq \frac{3}{2},\quad  {}^{t}(z,r),\ {}^{t}(z',r')\in \mathbb{R}^{2}_{+}.  \tag{5.2}
\end{align*}\\
In terms of $\zeta=\omega^{\theta}/r$, impulse and circulation (mass) can be written as

\begin{align*}
\frac{1}{2}\int_{\mathbb{R}^{3}}x\times \omega \dd x=\pi\left(\int_{\mathbb{R}^{2}_{+}}r^{3}\zeta\dd z\dd r\right)e_{z},\quad \int_{\{r=0\}}u\cdot \dd l(x)=\int_{\mathbb{R}^{2}_{+}}r\zeta \dd z\dd r.
\end{align*}\\
Besides them, we use a weighted $L^{1+\beta}$-norm to estimate the stream function.

\vspace{5pt}

%prop5.1
\begin{prop}
Let $0<\beta\leq 1$ and $0<\delta<1$. The estimate

\begin{align*}
|\psi(z,r)|\leq C\min\left\{r,\ \frac{1}{r^{1-\delta}}\right\}\left(||r^{3}\zeta||_1+||r\zeta||_1+||r^{1+2\beta}\zeta^{1+\beta}||_{1}^{1/(1+\beta)}\right)
,\ {}^{t}(z,r)\in \mathbb{R}^{2}_{+},   \tag{5.3}
\end{align*}\\
holds for $\psi$ in (5.1) with some constant $C$. 
\end{prop}

\vspace{5pt}

\begin{proof}
We set $s=\sqrt{|z-z'|^{2}+|r-r'|^{2} }$ and 

\begin{align*}
|\psi(z,r)|\leq \int_{\mathbb{R}^{2}_{+}}G(z,r,z',r')|\zeta(z',r')|r'\dd z'\dd r'
=\int_{s<r/2}+\int_{s\geq r/2}=: I+II.
\end{align*}\\
Since ${r'}^{2}/s^{3}\leq C/r$ for $s\geq r/2$, by the estimate of the Green function (5.2) for $\tau=3/2$, 

\begin{align*}
&II\lesssim \int_{s\geq r/2}\frac{(rr')^{2}}{s^{3}}|\zeta(z',r')|r'\dd z'\dd r'
\lesssim r||r\zeta||_{1},\\
&II \lesssim \frac{1}{r}\int_{s\geq r/2}r'^{2}|\zeta(z',r')|r'\dd z'\dd r'
\lesssim \frac{1}{r}||r^{3}\zeta||_{1}.
\end{align*}\\
Thus, $II\lesssim \min\{r,r^{-1}\}(||r\zeta||_{1}+ ||r^{3}\zeta||_{1})$. By (5.3) for $0<\tau\leq 3/2$ and the H\"older's inequality, for $0<\alpha\leq \beta\leq 1$,

\begin{align*}
I
&\lesssim r^{1/2+\tau}\int_{s<r/2}\frac{(r')^{1/2+\tau}}{s^{2\tau}}{r'}^{2}\zeta(z',r')\frac{1}{r'}\dd z'\dd r'\\
&\leq r^{1/2+\tau}\left(\int_{s<r/2}\frac{(r')^{(1/2+\tau)\sigma}}{s^{2\tau\sigma}}\frac{1}{r'}\dd z'\dd r'\right)^{1/\sigma}\left(\int_{s<r/2}({r'}^{2}\zeta(z',r'))^{1+\alpha}\frac{1}{r'}\dd z'\dd r'\right)^{1/(1+\alpha)}\\
&=Cr^{2-1/(1+\alpha)}||r^{2}\zeta||_{L^{1+\alpha}(\{s<r/2\}; r^{-1})  },
\end{align*}\\
with some constant $C$, independent of $r$. The constant $\sigma$ is the  H\"older conjugate to $1+\alpha$. We chose $\tau<\alpha/(1+\alpha)$ so that the integral is finite. By the H\"older's inequality,

\begin{align*}
||r^{2}\zeta||_{L^{1+\alpha}(\{s<r/2\}; r^{-1})  }\leq ||r^{2}\zeta||_{L^{1+\beta}(\{s<r/2\}; r^{-1})  }^{\theta}||r^{2}\zeta||_{L^{1}(\{s<r/2\}; r^{-1})  }^{1-\theta},\quad \frac{1}{1+\alpha}=\frac{\theta}{1+\beta}+1-\theta.
\end{align*}\\
Since $||r^{2}\zeta||_{L^{1+\beta}(\mathbb{R}^{2}_{+}; r^{-1})}=||r^{1+2\beta}\zeta^{1+\beta}||_{1}^{1/(1+\beta)}$ and $||r^{2}\zeta||_{L^{1}(\{s<r/2\}; r^{-1})  }
\lesssim \min\left\{1,r^{-2}\right\}(||r\zeta||_{1}+ ||r^{3}\zeta||_{1})$, we obtain

\begin{align*}
I\lesssim \min\left\{r^{2-1/(1+\alpha)  }, \frac{1}{r^{1/(1+\alpha)-2\theta }}   \right\}\left(||r^{3}\zeta||_1+||r\zeta||_1+||r^{1+2\beta}\zeta^{1+\beta}||_{1}^{1/(1+\beta)}\right).
\end{align*}\\
The right-hand side is $O(r^{2-1/(1+\beta)})$ as $r\to0$ for $\alpha=\beta$ and $O(r^{-1+\delta})$ as $r\to\infty$ for sufficiently small $\alpha$. Thus 

\begin{align*}
I\lesssim \min\left\{r^{2-1/(1+\beta)  }, \frac{1}{r^{1-\delta }}   \right\}\left(||r^{3}\zeta||_1+||r\zeta||_1+||r^{1+2\beta}\zeta^{1+\beta}||_{1}^{1/(1+\beta)}\right).
\end{align*} \\
By combining the estimates for $I$ and $II$, the desired estimate follows.
\end{proof}

\vspace{5pt}

%prop5.2
\begin{prop}
The estimate (5.3) holds for solutions to the Dirichlet problem

\begin{equation*}
\begin{aligned}
-L\psi&=r^{2}\zeta\quad \textrm{in}\ D(R), \\
\psi&=0\hspace{19pt} \textrm{on}\ \partial D(R),
\end{aligned}
\end{equation*}
with some constant independent of $R$.
\end{prop}

\vspace{5pt}

\begin{proof}
The function $\varphi=\psi/r^{2}$ is a solution to $-\Delta_y \varphi= \zeta$ in $B$ and $\varphi=0$ on $\partial B$. We may assume that $\zeta$ is non-negative. By zero extension of $\zeta$, we set $\tilde{\psi}$ by (5.1). Then $\tilde{\varphi}=\tilde{\psi}/r^{2}$ is positive and satisfies $-\Delta_y \tilde{\varphi}= \zeta$ in $B$. Since $\varphi-\tilde{\varphi}$ is harmonic in $B$ and negative on $\partial B$, $\varphi< \tilde{\varphi}$ by the maximum principle. Thus the result follows from Proposition 5.1.     
\end{proof}

\vspace{5pt}

We apply the estimate (5.3) to the grand state constructed in Lemma 4.9. The following energy identity is essentially due to Friendman-Turkington \cite[Lemma 3.2.]{FT81}.

\vspace{5pt}

%prop5.3
\begin{prop}
The identity 

\begin{align*}
\frac{1}{2\pi^{2}}||\psi||_{H}^{2}=\frac{1}{\lambda^{\beta}}||r^{1+2\beta}\zeta^{1+\beta}||_{1}+||r^{3}\zeta ||_{1}+\gamma ||r\zeta ||_{1}    \tag{5.4}
\end{align*}\\
holds for solutions to (3.5) with $\zeta=\lambda r^{-2}(\psi-r^{2}-\gamma)_{+}^{2q-1}$ and $\beta=(2q-1)^{-1}\in (0,1)$.
\end{prop}

\vspace{5pt}

\begin{proof}
By multiplying $\Psi=\psi-r^{2}-\gamma$ by $r\zeta $ and integrating it on $D$, we have 

\begin{align*}
\int_{D}\psi r \zeta \dd z\dd r=\int_{D}\Psi r\zeta \dd z\dd r+||r^{3}\zeta||_{1}+\gamma ||r\zeta||_1.
\end{align*}\\
Since $\zeta=\lambda r^{-2}\Psi_{+}^{2q-1}$, $\Psi_{+}=(\lambda^{-1}r^{2}\zeta)^{\beta}$ for $\beta=(2q-1)^{-1}$ and 

\begin{align*}
\int_{D}\Psi r\zeta\dd z\dd r
=\frac{1}{\lambda^{\beta}}\int_{D}r^{1+2\beta}\zeta^{1+\beta}\dd z\dd r
=\frac{1}{\lambda^{\beta}}||r^{1+2\beta}\zeta^{1+\beta}||_1.
\end{align*}\\
Since $\varphi=\psi/r^{2}$ satisfies $-\Delta_{y}\varphi=\zeta$ in $B$, $\varphi=0$ on $\partial B$ and $H(D;r^{-1})\cong F(B)$,  

\begin{align*}
\int_{D}\psi r \zeta \dd z\dd r=-\int_{D}\varphi\Delta_{y}\varphi r^{3}\dd z\dd r
=-\frac{1}{2\pi^{2}}\int_{B}\varphi\Delta_y\varphi\dd y
=\frac{1}{2\pi^{2}}\int_{B}|\nabla_y \varphi|^{2}\dd y=\frac{1}{2\pi^{2}}||\psi||_{H}^{2}.
\end{align*}\\
We obtained (5.4).
\end{proof}

\vspace{5pt}

%lem5.4
\begin{lem}
The grand sate in Lemma 4.9 satisfies

\begin{align*}
|\psi (z,r)|\leq C\min\left\{r,\ \frac{1}{r^{1-\delta}}\right\},\quad {}^{t}(z,r)\in D(R),   \tag{5.5}
\end{align*} \\
for $0<\delta<1$. The constant $C$ depends only on $||\psi||_{H}$, $\lambda$ and $\gamma$.
\end{lem}

\vspace{5pt}

\begin{proof}
The result follows from Propositions 5.2 and 5.3.
\end{proof}

%5.2
\subsection{Uniform boundedness of the vortex core}
We take an increasing sequence $\{R_n\}$ and set $H(D_n;r^{-1})$, $N(D_n;r^{-1})$ and $c_n=\inf_{N(D_n;r^{-1})}I$ with $D_n=D(R_n)$. By the zero extension, 

\begin{align*}
&H(D_n;r^{-1})\subset H(D_{n+1};r^{-1})\subset \cdots \subset H(\mathbb{R}^{2}_{+};r^{-1}),\\
&N(D_n;r^{-1})\subset N(D_{n+1};r^{-1})\subset \cdots \subset N(\mathbb{R}^{2}_{+};r^{-1}),\\
&c_n\geq c_{n+1}\geq \cdots \geq c.
\end{align*}

\vspace{5pt}

%prop5.9
\begin{prop}
\begin{align*}
\lim_{n\to\infty}c_n= c.    \tag{5.6}
\end{align*}
\end{prop}

\vspace{5pt}

\begin{proof}
We prove that for arbitrary $\phi\in N(\mathbb{R}^{2}_{+};r^{-1})$ there exists $\phi_n\in N(D_n;r^{-1})$ such that $\phi_n\to \phi$ in $H(\mathbb{R}^{2}_{+};r^{-1})$. This implies that 

\begin{align*}
c\leq \lim_{m\to\infty}c_n
=\lim_{n\to\infty}\inf_{N(D_n;r^{-1})}I
\leq \lim_{n\to\infty}I[\phi_n]=I[\phi].
\end{align*}\\
By taking infimum for $\phi$, (5.6) follows. 

We construct $\tilde{\phi}_{n}\in H(\mathbb{R}_{+}^{2}; r^{-1})$ supported in $D_n$ such that $\tilde{\phi}_n\to \phi$ in $H(\mathbb{R}^{2}_{+}; r^{-1})$ by a cut-off function argument. By Proposition 4.4, for $g_n(t)=I[t\tilde{\phi}_n]$ there exists $t(\tilde{\phi}_n)>0$ such that 

\begin{align*}
0=\frac{\dot{g}_{n}(t(\tilde{\phi}_n))}{t(\tilde{\phi}_n)}=||\tilde{\phi}_n||_{H(D_n;r^{-1})}^{2}-2\pi^{2}\lambda t(\tilde{\phi}_n)^{2(q-1)}\int_{D_n}\left(\tilde{\phi}_n-\frac{r^{2}+\gamma}{t(\tilde{\phi}_n)}\right)_{+}^{2q-1}\tilde{\phi}_n \frac{1}{r}\dd z\dd r.
\end{align*}\\
Since $\lim_{n\to\infty}||\tilde{\phi}_n||_{H(D_n;r^{-1})}=||\phi||_{H(\mathbb{R}^{2}_{+};r^{-1})}\neq 0$, the sequence $\{t(\tilde{\phi}_n)\}$ is bounded. 

Suppose that $\{t(\tilde{\phi}_n)\}$ does not converge to $1$. Then, there exists a subsequence such that $t(\tilde{\phi}_n)\to t_0$ for some $t_0> 0$. Sending $n\to\infty$ implies that $\dot{g}(t_0)=0$ for $g(t)=I[t\phi]$. Since Proposition 4.4 holds also for $\mathbb{R}^{2}_{+}$ and $\phi\in N(\mathbb{R}^{2}_{+};r^{-1})$, $t_0=t(\phi)=1$. We thus conclude that $t(\tilde{\phi}_n)\to 1$.

The desired sequence is obtained by setting $\phi_n=t(\tilde{\phi}_n)\tilde{\phi}_n\in N(D_n; r^{-1})$.
\end{proof}

\vspace{5pt}

The grand state $\psi_n$ of $c_n=\inf_{N(D_n;r^{-1})}I$ is uniformly bounded and equi-continuous in $\overline{\mathbb{R}^{2}_{+}}$ with a bounded vortex core. 

\vspace{5pt}

%prop5.1
\begin{prop}
The grand state $\psi_n\in N(D_n;r^{-1})$ in Lemma 4.9 satisfies  

\begin{align*}
\sup_{n}||\psi_n||_{H(D_n;r^{-1})}<\infty.    \tag{5.7}
\end{align*}
\end{prop}

\vspace{5pt}

\begin{proof}
By (4.2),

\begin{align*}
\frac{1}{2}\left(1-\frac{1}{q}\right)||\psi_n||_{H(D_n;r^{-1})}^{2}
\leq I[\psi_n]=c_n\leq c_1.
\end{align*}
\end{proof}

%prop5.7
\begin{prop}
There exists some $\nu\in (0,1)$ such that $\varphi_n=\psi_n/r^{2}$ satisfies 

\begin{align*}
\sup_{n}||\varphi_n||_{C^{\nu}(\mathbb{R}^{5})}<\infty.   \tag{5.8}
\end{align*}
\end{prop}

\vspace{5pt}

\begin{proof}
By (5.5) and (5.7), 

\begin{align*}
|\psi_n(z,r)|\leq C\min\left\{r,\ \frac{1}{r^{1-\delta}}\right\},\quad {}^{t}(z,r)\in D_n,
\end{align*} \\
for $0<\delta<1$, with some constant $C$, independent of $n$. Thus $\zeta_n=\lambda r^{-2}(\psi_n-r^{2}-\gamma)_+^{2q-1}$ satisfies $\zeta_n\leq C(1+r^{-1})$, $(z,r)\in D_n$. Hence, 

\begin{align*}
\sup_{n}||\zeta_n||_{L^{s}_{\textrm{ul}}(\mathbb{R}^{5}) }<\infty,\quad \frac{5}{2}<s<4,
\end{align*}\\
where $L^{s}_{\textrm{ul}}(\mathbb{R}^{5})$ denotes the uniformly local $L^{s}$ space on $\mathbb{R}^{5}$. Since $-\Delta \varphi_n=\zeta_n$ in $B_n$ and $\varphi_n=0$ on $\partial B_n$, by the elliptic regularity, uniformly local $L^{s}$-norm of $\varphi_n$ is uniformly bounded up to second orders. Thus (5.8) holds for $\nu=1-5/s^{*}$ and $1/s^{*}=1/s-1/5$ by the Sobolev embedding. 
\end{proof}

\vspace{5pt}

The equi-continuity of $\{\varphi_n\}$ implies that the vortex core is uniformly bounded \cite[Lemma 4.1]{AS89}.

\vspace{5pt}

%prop5.5
\begin{prop}
There exists $R>0$ such that $\Omega_n=\{{}^{t}(z,r)\in \mathbb{R}^{2}_{+}\ |\ \psi_n-r^{2}-\gamma>0\}\subset D(R)$ for all $n\geq 1$.
\end{prop}

\vspace{5pt}

\begin{proof}
For notational simplicity, we write $\varphi_n(z,r)=\varphi_n(y)$. By the Sobolev inequality $F(\mathbb{R}^{5})\subset L^{10/3}(\mathbb{R}^{5})$, 

\begin{align*}
|\{z\in \mathbb{R}\ |\ \varphi_n(z,r_0)\geq 1/2 \ \}|
&\lesssim \int_{-\infty}^{\infty}\varphi_n^{8/3}(z, r_0)\dd z\\
&=-\int_{r_0}^{\infty}\frac{\dd }{\dd r}\int_{-\infty}^{\infty}\varphi_n^{8/3}\dd z\dd r\\
&\lesssim \frac{1}{r_0^{3}}\int_{0}^{\infty}\int_{-\infty}^{\infty}|\nabla \varphi_n|\varphi_n^{5/3}r^{3}\dd z\dd r
\lesssim \frac{1}{r_0^{3}}||\varphi_n||_{F(B_n)}^{8/3}.
\end{align*}\\
Thus by (5.7) and $H(D_n;r^{-1})\cong F(B_n)$, 

\begin{align*}
\left|\left\{z\in \mathbb{R}\ \middle|\ \varphi_n(z,r_0)\geq 1/2\ \right\}\right|\leq \frac{C}{r^{3}_0},\quad r_0>0,   
\end{align*}\\
with some constant $C$, independent of $n$. We take $r_n>0$ such that $r_n=\sup\ \{r\ |\ {}^{t}(z,r)\in \Omega_n\ \}$. The maximum point lies on the $r$-axis and $\partial\Omega_n$ since $\Omega_n$ is symmetric for $z$. Thus $\varphi_n(0,r_n)=1+r_n^{-2}\gamma\geq 1$. Since $\{\varphi_n\}$ is equi-continuous around ${}^{t}(0,r_n)$ by (5.8), there exists $\delta_1>0$ such that 

\begin{align*}
\varphi_n(z, r_n)\geq \frac{1}{2},\quad |z|\leq \delta_1. 
\end{align*}\\
Thus $r_n\lesssim \delta^{-1/3}_1$ and $\Omega_n$ is bounded in the $r$-direction. 

In a similar way, we take $z_n>0$ such that $z_n=\sup\ \{z\ |\ {}^{t}(z,r)\in \Omega_n\ \}$ and the maximum point ${}^{t}(z_n,r_n)$. Since $\varphi_n(z_n, r_n)\geq 1$, by equi-continuity of $\varphi_n$, there exists $\delta_2>0$ such that 

\begin{align*}
\varphi_n(z_n, r_n+\delta_2)\geq \frac{1}{2}.
\end{align*}\\
Since $\varphi_n$ is decreasing for $z$, $z_n\lesssim \delta_2^{-3}$ and $\Omega_n$ is bounded in the $z$-direction.
\end{proof}

\subsection{Convergence to a grand state in $\mathbb{R}^{2}_{+}$}

%lem5.8
\begin{lem}
The sequence $\{\varphi_n\}$ subsequently converges to a limit $\varphi$ locally uniformly in $\mathbb{R}^{5}$. The limit $\psi=r^{2}\varphi\in N(\mathbb{R}^{2}_{+};r^{-1})$ is a grand state to (3.4) such that $I'[\psi]=0$ and $\psi=\psi^{*}$ with a bounded vortex core $\Omega=\{{}^{t}(z,r)\in \mathbb{R}^{2}_{+}\ |\ \psi-r^{2}-\gamma>0\}$. 
\end{lem}

\vspace{5pt}

\begin{proof}
The sequence $\{\varphi_n\}$ is uniformly bounded and equi-continuous by the uniform estimate (5.8). By choosing a subsequence, $\varphi_n$ converges to a limit $\varphi$ locally uniformly in $\mathbb{R}^{5}$ by Ascoli-Arzel\'a theorem. By Proposition 5.8, $U_n=\{y={}^{t}(y_1,y')\in \mathbb{R}^{5}\ |\ y_1=z,\ |y'|=r,\ {}^{t}(z,r)\in \Omega_n \}$ is uniformly bounded. Thus, $\zeta_n=\lambda r^{4(q-1)}(\varphi_n-1-r^{-2}\gamma)_{+}^{2q-1}$ converges to $\zeta=\lambda r^{4(q-1)}(\varphi-1-r^{-2}\gamma)_{+}^{2q-1}=\lambda r^{-2}(\psi-r^{2}-\gamma )_{+}^{2q-1}$ uniformly in $\mathbb{R}^{5}$ and $\Omega=\{ \psi-r^{2}-\gamma>0\}$ is bounded. By $\textrm{spt}\ \omega^{\theta}=\textrm{spt}\ \zeta=\overline{\Omega}$, $\Omega$ is the vortex core. By the uniform estimate of the Dirichlet energy (5.7), $\psi\in H(\mathbb{R}^{2}_{+};r^{-1})$. The properties $I'[\psi]=0$ and $\psi=\psi^{*}$ follow from those of $\psi_n$. 

We show that the limit $\varphi$ is non-trivial. Suppose that $\varphi\equiv 0$. Then, $0\leq \varphi_n\leq 1$ in $B_n$ for large $n$. By Lemma 3.2, $\varphi_n$ solves $-\Delta_y \varphi_n=\lambda r^{4(q-1)}(\varphi_n-1-r^{-2}\gamma )_{+}^{2q-1}=0$ in $B_n$ and $\varphi_n=0$ on $\partial B_n$. Thus $\varphi_n\equiv 0$. This contradicts $\psi_n\nequiv 0$. Thus the limit is non-trivial. In particular, $\psi\in N(\mathbb{R}^{2}_{+};r^{-1})$.

It remains to show that $\psi$ is a grand state to (3.4). By the uniform boundedness of the vortex core $\Omega_n=\{\psi_n-r^{2}-\gamma>0\}\subset D(R)$, 

\begin{align*}
J'[\psi_n]\psi_n
&=2\pi^{2}\lambda\int_{D(R)}(\psi_n-r^{2}-\gamma)_{+}^{2q-1}\psi_n\frac{1}{r}\dd z\dd r \\
&\to 2\pi^{2}\lambda\int_{D(R)}(\psi-r^{2}-\gamma)_{+}^{2q-1}\psi\frac{1}{r}\dd z\dd r
=J'[\psi]\psi \quad \textrm{as}\ n\to\infty.
\end{align*}\\
By $0=I'[\psi_n]\psi_n=||\psi_n||_{H(\mathbb{R}^{2}_{+};r^{-1})}^{2}+J'[\psi_n]\psi_n$ and $0=I'[\psi]\psi=||\psi||_{H(\mathbb{R}^{2}_{+};r^{-1})}^{2}+J'[\psi]\psi$, we have $\lim_{n\to\infty}||\psi_n||_{H(\mathbb{R}^{2}_{+};r^{-1} )}=||\psi||_{H(\mathbb{R}^{2}_{+};r^{-1} )}$. By choosing a subsequence, $\psi_n\to \psi$ in $H(\mathbb{R}^{2}_{+};r^{-1} )$. By continuity of $I$ and (5.6), 

\begin{align*}
\inf_{N(\mathbb{R}^{2}_{+};r^{-1})}I=c=\lim_{n\to\infty}c_n=\lim_{n\to\infty}I[\psi_n]=I[\psi].
\end{align*}\\
Thus $\psi$ is a grand state to (3.4). 
\end{proof}

\vspace{5pt}

\begin{proof}[Proof of Theorem 2.1]
The critical point $\psi\in H(\mathbb{R}^{2}_{+};r^{-1})$ is a classical solution $\psi\in C^{2+\nu}(\overline{\mathbb{R}^{2}_{+}})$ of (2.5) for $\nu\in (0,1)$ and $\psi$ is decreasing for $|z|$ by a maximum principle as we proved Lemma 3.2 and Proposition 4.10 for a grand state in a half disk. The level sets $\{(\psi-r^{2}-\gamma)_{+}=l\}$ are closed curves of class $C^{2+\nu}$ for $0\leq l<l_0$ and $l_0=||(\psi-r^{2}-\gamma)_{+}||_{\infty}$ and points for $l=l_0$. The vortex core $\Omega$ is bounded, connected and simply-connected as we proved Lemma 4.12 for a grand state in $D$. Since $r^{-1}\nabla_{z,r}G(z,r,z',r')\to0$ as $z^{2}+r^{2}\to\infty$ for each ${}^{t}(z',r')\in \mathbb{R}^{2}_{+}$ and $\zeta=\lambda r^{-2}(\psi-r^{2}-\gamma)_{+}^{2q-1}$ is supported in $\overline{\Omega}$, $r^{-1}\nabla_{z,r}\psi\to 0$ follows from (5.1).
\end{proof}

\vspace{5pt}

\begin{proof}[Proof of Theorem 1.1]
By rotational invariance of (1.4), we may assume that $u_{\infty}={}^{t}(0,0,-W)$ for $W>0$. We take $1<q<\infty$ and $0<\lambda,\gamma<\infty$. Then by Theorem 2.1, there exists a solution $\psi\in C^{2+\nu}(\overline{\mathbb{R}^{2}_{+}})$ of (2.5) for $\nu\in (0,1)$. We set $f=\dot{\Gamma}(\Psi)$ by $\Psi=\psi-Wr^{2}/2-\gamma$ and $\dot{\Gamma}(t)=(\lambda q)^{1/2}t_{+}^{q-1}$. Then, $f\in C(\overline{\mathbb{R}^{2}_{+}})$. By regarding it to an axisymmetric function in $\mathbb{R}^{3}$, $f\in C(\mathbb{R}^{3})$. If $q>3$, $f\in C^{2+\mu}(\mathbb{R}^{3})$ for some $\mu\in (0,1)$. 

The velocity $u$ defined by $(2.1)$ satisfies $u\in C^{1}(\mathbb{R}^{3}\backslash \{r=0\})$ and $\nabla \times u=f u$ by (2.2). Since $f$ is supported on a solid torus rotating $\overline{\Omega}$ around the $z$-axis, $u$ is harmonic near the $z$-axis. By (2.1) and (2.5),  

\begin{align*}
u&=-\frac{1}{r}\partial_z\Psi e_{r}(\theta)+\frac{1}{r}\Gamma(\Psi)e_{\theta}(\theta)+\frac{1}{r}\partial_r\Psi e_{z} \\
&=-\frac{1}{r}\partial_z\psi e_{r}(\theta)+\frac{1}{r}\Gamma(\Psi)e_{\theta}(\theta)+\frac{1}{r}\partial_r\psi e_{z}+u_{\infty}\to u_{\infty}\quad \textrm{as}\ |x|\to\infty. 
\end{align*}\\
Thus $u\in C^{1}(\mathbb{R}^{3})$ is an axisymmetric solution to (1.4). Since the level sets $\{(\psi-Wr^{2}/2-\gamma)_{+}=l \}$ for $0<l<l_0$ are symmetric nested closed curves in $\mathbb{R}^{2}_{+}$ (resp. points for $l=l_0$), the level sets $f^{-1}(k)$ are symmetric nested tori for $0<k<k_0$ with $k_0=\sqrt{\lambda q}l_0^{q-1}$ (resp. circles for $k=k_0$). The torus $f^{-1}(k)=\Sigma(l)$ includes vortex lines (stream lines) associated with $l=(k/\sqrt{\lambda q})^{1/(q-1)}$ for $0<k<k_0$. All vortex lines are closed or quasi-periodic if $\Theta(l)$ is commensurable with $2\pi$ or not. Thus $f^{-1}(k)$ is an invariant torus. The proof is now complete.
\end{proof}

\section*{Acknowledgements}
This work is partially supported by JSPS through the Grant-in-aid for Young Scientist 20K14347, Scientific Research (B) 17H02853 and MEXT Promotion of Distinctive Joint Research Center Program Grant Number JPMXP0619217849.\\

%ref
\bibliographystyle{plain}
\bibliography{ref}

\end{document}